\documentclass[12pt]{article}
\usepackage{amsfonts}
\usepackage{epsfig}
\title{The $5$-Electron case of Thomson's Problem}
\author{Richard Evan Schwartz \thanks{\hskip 5 pt Supported by 
N.S.F. Research Grant DMS-0072607}}

\newtheorem{theorem}{Theorem}[section]

\newtheorem{lemma}[theorem]{Lemma}

\newtheorem{corollary}[theorem]{Corollary}

\def\startproof{{\bf {\medskip}{\noindent}Proof: }}

\def\endproof{$\spadesuit$  \newline}

\def\C{\mbox{\boldmath{$C$}}}%
\def\D{\mbox{\boldmath{$D$}}}%
\def\Q{\mbox{\boldmath{$Q$}}}%
\def\R{\mbox{\boldmath{$R$}}}%

\begin{document}
\maketitle
\begin{abstract}
We give a rigorous, computer-assisted proof that the
triangular bi-pyramid is the unique configuration of
$5$ points on the sphere that globally 
minimizes the Coulomb $(1/r)$ potential.  We also
prove the same result for the $(1/r^2)$ potential.
The main mathematical contribution of the paper
is a fairly efficient energy estimate that works for
any number of points and any power law potential.
\end{abstract}

\section{Introduction}

\subsection{Background}

The problem of finding how electrons optimally distribute themselves
on the sphere is a well-known and difficult one.  It is known as
{\it Thomson's problem\/}, and dates from J. J. Thomson's $1904$
publication [{\bf T\/}].    Thomson's problem is of interest not just
to mathematicians, but also to physicists and chemists.

Here is a mathematical
formulation.
Let $S^2 \subset \R^3$ be the unit
sphere.   Let $P$ be a collection of
$n$ distinct point $p_1,...,p_n \in S^2$.
Let $E: \R^+ \to \R^+$ be some function.
the total energy to be the sum
\begin{equation}
{\cal E\/}(P)=\sum_{i>j} E(\|p_i-p_j\|).
\end{equation}
Here $\|\cdot\|$ is the usual norm on $\R^3$, and
$\|p_i-p_j\|$ is the distance from $p_i$ to $p_j$.
The question is then:  What does $P$ look like if ${\cal E\/}(P)$ is
as small as possible?  

The question is perhaps too broad as stated, because the answer
likely depends on the function $E$.    To consider a narrower question,
one restricts the class of functions in some way.  For instance,
a natural class of potential functions is given by the
{\it power laws\/}
\begin{equation}
E(r)=\frac{1}{r^e}; \hskip 30 pt e \in (0,\infty).
\end{equation}
The case $e=1$ is specially interesting to physicists.  It is known as
the {\it Coulomb potential\/}.  When energy is measured
with respect to the Coulomb potential, the points are naturally
considered to be electrons.  

There is a large literature
on Thomson's problem.  One early work on Thomson's problem is
[{\bf C\/}].
The paper [{\bf SK\/}] gives a nice survey in the two dimensional case,
with an emphasis on the case when $n$ is large. 
The paper [{\bf BBCGKS\/}] gives a survey
of results, both theoretical and experimental, about
highly symmetric configurations in higher dimensions.
The fairly recent paper [{\bf RZS\/}] has
some theoretical bounds for the logarithmic potential,
and also has a large amount of experimental information
about configurations minimizing the power laws on
the $2$-sphere. 
 The website [{\bf CCD\/}] has a
list of experimentally determined (candidate) minimizers for
the Coulomb potential for $n=2,...,972$.   

There are certain values of $n$ where the  minimal configuration
is rigorously known for all the power laws.  
\begin{itemize}
\item When $n=2$, the points of $P$ are antipodal.
\item When $n=3$, the points of $P$ make an 
equilateral triangle in an equator.
\item  When
$n=4$, the points of $P$ are the vertices of a regular
tetrahedron. 
\item When $n=6$, the points of $P$ are the vertices of a regular
octahedron.
\item 
 When $n=12$, the points of $P$ are the
vertices of a regular icosahedron.
\end{itemize}
The cases $n=2$ is trivial to prove, and the case $n=3$ is an
easy exercise.  The cases
$n=4,6,12$ are all covered in [{\bf CK\/}, Theorem 1.2],
a much broader result concerning a fairly general
kind of energy potential and a class of
point configurations (not restricted to $2$ dimensions) called
{\it sharp configurations\/}.  For the specific cases
cases $n=4,6$ see the older work [{\bf Y\/}].

The case $n=5$ is conspicuously absent from the above list of
known results. 
Everyone agrees, based on numerical simulation,
 that the triangular bi-pyramid (TBP)
 seems to be the global minimizer
for the Coulomb potential.
In the TBP, two points are antipodal points on $S^2$  and
the remaining $3$ points form an equilateral
triangle on the equator midway between
the two antipodal points. 
More generally, numerical experiments \footnote{Lacking a handy reference
for these experiments, we performed our own.} suggest the following.
\begin{itemize}
\item The TBP is a local minimizer for the power law potential with exponent $e$ if
and only if $e \in (0,e_1)$, with $e_1 \approx 21.147123$.
\item The TBP is a global minimizer for the power law potential with exponent $e$
if and only if $e \in (0,e_2)$, with $e_2 \approx 15.040808$.
\end{itemize}
For $e>e_2$, it seems that the global minimizer is a pryamid with
square base.  The precise configuration depends on $e$ in this range.

In spite of detailed experimental knowledge about the case $n=5$, 
it seems there has not ever
been a proof that the TBP is
global minimum for any power law potential.
In particular, this has not been proved
for the Coulomb potential.
As far as we know, there are two rigorous results
for the case $n=5$.
\begin{itemize}
\item The paper [{\bf DLT\/}] contains a (traditional) proof
that the TBP maximizes the geometric mean of the pairwise
distances between the points.  This case corresponds to the
logarithmic potential $E(r)=-\log(r)$.
\item The paper [{\bf HS\/}]
contains a computer-aided proof that
the TBP {\it maximizes\/}
the potential for the exponent $e=-1$.   That is, $5$ points on the sphere
arrange themselves into a TBP so as to maximize the
total sum of the pairwise distances.   
\end{itemize}

\subsection{Results}

It is the purpose of this paper to prove the following result.
\begin{theorem}
\label{main}
The TBP is the unique configuration of minimal
energy with respect to the Coulomb  potential.
\end{theorem}

Our proof is computer aided, and similar in spirit to [{\bf HS\/}].  
While the argument in [{\bf HS\/}] is exactly tailored to
understanding the sums of the distances, our method
is rather insensitive to the
precise power law being used.  Just to illustrate this fact, we prove

\begin{theorem}
\label{main2}
The TBP is the unique configuration of minimal
energy with respect to the $1/r^2$ potential.
\end{theorem}

We could certainly add other exponents to our list of
results.  However,  we currently have implemented the
interval arithmetic in such a way that it works just for
the exponents $e=1,2$.  We did this to avoid using
the pow function ${\rm pow\/}(a,b)=a^b$, which
is not covered by the same guarantees in the
IEEE standards as are the basic arithmetic
operations.  See [{\bf I\/}] and [{\bf I2\/}].
Mainly what stops us is a sense of diminishing returns.

However, with a view towards an eventually broader application,
we try to state as many estimates as we can for the
general power law.   In every situation, we try to
point out the exact
generality with which the construction holds.
For example, our main technical result,
Theorem \ref{EE}, holds for any function $E$
satisfies the conditions discussed in \S \ref{conditions}.
Also, the version of
Theorem \ref{EE} we prove works for general
$n$-point configurations.

\subsection{Outline of the Proof}

It has probably been clear for a long time that one
could prove a result like Theorem \ref{main} using a
computer program.  The main difficulties are
technical rather than conceptual.  The hard part is
getting estimates that are sharp enough
so as to lead to a feasible calculation. 
Our paper doesn't really have any dramatic new ideas.
It just organizes things well enough to get the job done.

We begin by eliminating some obviously bad configurations
from consideration.  For example we show that no
two points in a minimizing configuration lie within 
$1/2$ units of each other.  This result holds for any
power law potential.    See Lemma \ref{coulomb}.
After eliminating these bad configurations, we are left with
a compact configuration space $\Omega$ of possible minimizers.

We use stereographic projection to transfer the 
configurations on $S^2$ to configurations in
$\C \cup \infty$.  By keeping one point at $\infty$
and another on the positive real axis, we end up
with a natural description of $\Omega$ as the
set $[0,4] \times [-2,2]^6 \subset \R^7$.
Working in $\Omega$, we define a natural way to
subdivide a rectangular solid subset of $\Omega$ into
smaller rectangular solids.   

We use a divide-and-conquer algorithm to show that any
minimizing configuration in $\Omega$ must lie very close
to the TBP.   Assuming that $\cal Q$ is a rectangular
solid that \underline{does not} lie too close to the
point(s) of $\Omega$ representing the TBP, we try to
eliminate $\cal Q$ with a $3$-step procedure.

\begin{itemize}

\item We eliminate $\cal Q$ if we can see, based on the
fact that the regular tetrahedron minimizes energy for
$4$ points, that no configuration in $\cal Q$ can
minimize energy.   This is described in \S \ref{tetra4}.
This method of eliminating $\cal Q$ cuts away a lot
of the junk, so to speak, and focuses our attention on
the configurations where are fairly near the TBP.

\item We eliminate $\cal Q$ if the configurations in it are
redundant, in the sense that some permutation of the vertices
or obvious application of symmetry changes these configurations
to ones in a form that we deem more standard.  See \S \ref{redundant0}.

\item If the first two methods fail, we
evaluate the energy $\cal E$ on the vertices of $\cal Q$ and then
apply an {\it a priori\/} estimate on how far $\cal E$ differs from
a linear function on $\cal Q$.    Our main result along these lines
is Theorem \ref{EE}.   This is our most powerful and general
method of elimination, but we only use it when all else fails.

\end{itemize}

If it is not possible to eliminate $\cal Q$ by any of our
methods, we subdivide $\cal Q$ into smaller
rectangular solids and try again on each piece.
This algorithm is discussed,
in outline, in \S \ref{discuss1}.   In the chapters following
\S \ref{discuss1}, we fill in the details.
Theorem \ref{EE} is the most
subtle and important part of the paper.  We
prove this result in \S 5-8.

The end result of our finite calculation is that the true minimizer of
$\cal Q$ must lie very close to the TBP.   See Lemma \ref{maincomp}.
  To finish the proof, we 
use calculus to show that the TPB can be the only local minimum
in the region where we have confined the minimizer.  We do this
by showing that the Hessian of $\cal E$, the matrix of
second partials, is positive definite throughout the region of
interest to us.  See Lemma \ref{mainhess}.

\subsection{Computational Issues}
\label{issues}

We implement our computer program in Java, using 
interval arithmetic to control the round-off errors. 
We will discuss this in \S \ref{interval}.
Our code only uses the basic operations
plus, minus, times, divide, and sqrt, and these operations
are performed in such a way as to avoid and arithmetic
errors (such as taking the square root of a negative number.)
It takes about $6$ hours for our program to eliminate all configurations
that are not within $2^{-14}$, in the
$L^{\infty}$ norm, of 
a suitably normalized version of the TBP.
See Lemma \ref{maincomp}.
Without the interval arithmetic, the code runs in about an hour.

Theorem \ref{EE} is the rate-determining step in our calculations.
Numerical experimentation suggests that our bound in Theorem \ref{EE}
is off by a factor somewhere between $2$ and $4$.    In light of this fact,
our calculation probably ought to be about $10$ times faster than it is.
We can certainly improve Theorem \ref{EE}, but we haven't been
been able to think of a simple or dramatic improvement. 

 A nice feature of our program is that we have
embedded it in a graphical user interface. The
reader can watch the program in action and see
how it samples the configuration space.  
the reader can also manually construct a rectangular solid subset
of the configuration and then see a printout
of all the computational tests that are applied to it.
This graphical aspect doesn't add anything to the
formal proof, but it makes it less likely we have
made a gross computing error.     We have tried
to isolate the relatively small amount of
computer code that goes into the actual proof,
so that it can be more easily inspected.

The entire Java program is available from my website.  See \newline
{\bf http://www.math.brown.edu/$\sim$res/Electron/index.html\/}. The code
is fairly well documented, and we're still working to improve the
documentation.  Aside from graphical support files,
all the files involved in the proof
have the {\bf Interval\/} prefix.  The directory contains a number
of other files, which support other features of the program.
Even though the proof portion of the code is done and working,
the whole program is still somewhat in flux.  I plan to gradually
improve the code as time passes, and update it as I go.

\subsection{Acknowledgements}

I would like to thank Henry Cohn for helpful conversations about
placing electrons on a sphere.  Henry's great colloquium talk at
Brown university this fall inspired me to work on this problem.
I would also like to thank Jeff Hoffstein and Jill Pipher for their
interest and encouragement while I worked on this problem.
I would like to thank John Hughes for a very interesting
discussion about interval arithmetic.  
Finally, I would like
to say that I learned how to do interval arithmetic
in Java by reading the source files from Tim Hickey's implementation,
{\bf http://interval.sourceforge.net/interval/index.html\/}.  

\newpage

\section{Proof in Broad Strokes}

\subsection{Stereographic Projection}

We find it convenient to work mainly with
$\C \cup \infty$ rather than on $S^2$.    
Our reason for this is that a configuration space
based on points in $\C$ has a
natural flat structure, and lends itself well
to a nice subdivision scheme.   The subdivision
scheme, which essentially amounts to cutting
rectangular solids into smaller rectangular
solids, feeds into our divide-and-conquer
algorithm.  All this is discussed in \S \ref{discuss1} below.

We map $S^2$ to $\C \cup \infty$ using
{\it stereographic projection\/}:
\begin{equation}
\label{stereo}
\Sigma(x,y,z)=\frac{x}{1-z}+ i\ \frac{y}{1-z}.
\end{equation}
$\Sigma$ is a conformal diffeomorphism which maps
circles on $S^2$ to generalized circles in $\C \cup \infty$.
A {\it generalized circle\/} is either a circle or a straight line.
We have 
\begin{equation}
\Sigma(0,0,1)=\infty.
\end{equation}
   Thus, $\Sigma$ maps
a circle $C \subset S^2$ to a straight line if and only if
$(0,0,1) \in C$.

The inverse map is given by
\begin{equation}
\label{INV}
\Sigma^{-1}(x+iy)=\bigg(\frac{2x}{1+x^2+y^2},\frac{2y}{1+x^2+y^2},
1-\frac{2}{1+x^2+y^2}\bigg).
\end{equation}

Let $\|\cdot\|$ denote the usual norm on $\R^3$.
Here are two pieces of metric information we will use later on:
\begin{equation}
\label{metric1}
\|\Sigma^{-1}(z)-\Sigma^{-1}(\infty)\|=\frac{2}{\sqrt{1+|z|^2}}
\end{equation}

\begin{equation}
\label{metric2}
\bigg{\|} \frac{d\Sigma^{-1}}{dx}\bigg{\|} =
\bigg{\|} \frac{d\Sigma^{-1}}{dy}\bigg{\|} =
\frac{2}{1+x^2+y^2}.
\end{equation} 

Equations \ref{metric1} and \ref{metric2} both have
straightforward derivations, which we omit.
Slightly abusing notation, we define
\begin{equation}
E(z_1,z_2)=E\big(\Sigma^{-1}(z_1),\Sigma^{-1}(z_2)\big)
\end{equation}
Here $E$ can be any energy potential.
for points $z_1,z_2 \in \C \cup \infty$.  In this way, we
can talk about the energy of a configuration of points in $\C \cup \infty$.

\subsection{The Triangular Bi-Pyramid}

As in the introduction $\|\cdot\|$ denotes the
usual norm on $\R^3$.  Unless stated otherwise,
all the distances we measure in $\R^3$ are taken
with respect to the Euclidean metric.  What we say here works for
any energy potential.

In the TBP, we have the following information.
\begin{itemize}
\item One pair of points is $2$ units apart.
\item $3$ pairs of points are $\sqrt 3$ units apart.
\item $6$ pairs of points are $\sqrt 2$ units apart.
\end{itemize}
Accordingly, the energy of the TBP, with respect to $E$, is
\begin{equation}
\label{penta}
M_E=E(1/2)+3E(\sqrt 3) + 6 E(\sqrt 2).
\end{equation}

It is well known that the regular tetrahedron minimizes the
energy for $4$ points.   All points are $\sqrt{8/3}$ units apart.
For this reason, the regular tetrahedron has energy
\begin{equation}
T_E= 6 E(\sqrt{8/3}).
\end{equation}
We can use this fact to get some crude bounds on
$5$-point configurations.   Given a $5$ point
configuration $P$ and a point $p \in P$, we define

\begin{equation}
\label{tetra1}
E(P,p)=\sum_{q \in P \{p\}} E(\|p-q\|).
\end{equation}
We have the immediate estimate
\begin{equation}
\label{tetra}
E(P) \geq E(P,p) + T_E.
\end{equation}
Thus, if $P$ is an energy minimizer we must have
\begin{equation}
\label{tetra3}
E(P,p) \leq M_E-T_E; \hskip 30 pt \forall p \in P.
\end{equation}
This is one of the criteria we will use to eliminate certain
configurations from consideration.  

We can also use Equation \ref{tetra3} to give information
about pairs of points within a minimizing configuration.
We will consider this in the next section.

\subsection{Estimates for the Power Laws}

The fairly weak results in this section are designed to give us some
control on the size of the configuration space we must consider.
For any given exponent, these results are just short calculations.
It is only our desire to handle all exponents at the same time that
adds complexity to the proof.

\begin{lemma}
\label{coulomb}
Let $p,q_1,q_2$ be $3$ points of an energy minimizer with
respect to a power law potential.  
Then $\|p-q_1\| > 1/2$.
\end{lemma}

\startproof
Let $E(r)=1/r^e$.  Since $E$ is decreasing,
$E(\|p-r\|) \geq E(2)$ for all $p,r\in S^2$,
because $S^2$ has chordal diameter $2$.   In light of
Equation \ref{tetra}, this lemma is true provided that
\begin{equation}
T_E+3E(2)+E(1/2)-M_E>0
\end{equation}
When $E$ is as above, our problem boils down to showing that
\begin{equation}
\phi(e)=2^{1-e}-3^{1-e/2} + 2^e - 3^{1-e/2} + 2^{1-(3e)/2} 3^{1+e/2}>0.
\end{equation}
This is an exercise in calculus.  We compute
\begin{equation}
\phi(0)=0; \hskip 30 pt \phi'(0)>0; \hskip 30 pt
\phi''(0)>1.
\end{equation}
We compute $\phi'''(e)=A(e)-B(e)$, where
$$
A(e)=3 \log(2)^2 2^{-2-e/2} + 2^e \log(2)^3 + \frac{\log(3)}{8} 3^{1-e/2}>0;
$$
\begin{equation}
B(e)=2^{1-e}\log(2)^3+2^{-2-(3e)/2}3^{1+e/2}\log(8/3)^3<2.
\end{equation}
In short,
$\phi'''(e)>-2$.
Taylor's theorem with remainder now tells us that
\begin{equation}
\phi(e)>\frac{e^2}{2}-\frac{e^3}{3}.
\end{equation}
Hence $\phi(e)>0$ for $e \in (0,3/2)$.
A similar computation for $\phi'$ shows that 
\begin{equation}
\phi'(e)>-10.
\end{equation}
We compute  that
\begin{equation}
\phi(3/2+j/10)>1; \hskip 30 pt j=0,...,85.
\end{equation}
Combining the last two equations, we see that $\phi(e)>0$ for all
$e \in [3/2,10]$.    Finally, for $e>10$, the result is obvious.
\endproof

\begin{lemma}
\label{coulomb2}
Let $p,q_1,q_2$ be $3$ points of an energy minimizer with
respect to a power law potential.   Assume that
$E(\|p-q_2\|) \leq E(\|p-q_1\|)$.  Then
$\|p-q_1\| > 1/2$ and $\|p-q_2\| > 1$.
\end{lemma}

\startproof
We have the inequality
\begin{equation}
E(\|p-q_1\|)+E(\|p-q_2\|) \leq M_E-T_E-2E(2).
\end{equation}
Hence
\begin{equation}
E(\|p-q_2\|) \leq \frac{M_E-T_E-2E(2)}{2}.
\end{equation}
Establishing this inequality boils down to showing that
\begin{equation}
\phi(e)=2+2^{1-e}-(3)(2^{1-e/2})-2^{-e}-3^{1-e/2}+2^{1-(3e)/2}3^{1+e/2}>0.
\end{equation}
This time we compute
\begin{equation}
\phi(0)=0; \hskip 30 pt
\phi'(0)>0; \hskip 30 pt
\phi''(0)>1/4.
\end{equation}
Examining the terms of $\phi'''$, as in the previous lemma, we find that
\begin{equation}
\phi'''(e)>-12.
\end{equation}
Taylor's theorem now tells us that $\phi>0$ on $[0,1/16)$.
A similar computation shows
\begin{equation}
\phi'(e)>-10.
\end{equation}
We now compute that
\begin{equation}
\phi(1/16 + j/2000)>1/200; \hskip 30 pt j=0,...,1875.
\end{equation}
Combining the last two equations, we see that $\phi(e)>0$ for
$e \in [1/16,1]$.
Next, we compute that
\begin{equation}
\phi(1+ j/100)>1/10; \hskip 30 pt j=0,...,900.
\end{equation}
This shows that $\phi(e)>0$ for $e \in [1,10]$.  For
$e>10$ the result is again obvious.

\subsection{Planar Configurations}
\label{twomin}
\label{space}

We took the trouble to prove Lemmas \ref{coulomb} and
\ref{coulomb2} for any power law so that what we say
in this section works for any power law.
\newline
\newline
\noindent
{\bf Basic Definition:\/}
Let $P=\{p_0,...,p_4\}$ be a configuration of $5$ points on $S^2$.
$Z=\Sigma(P)=\{z_0,...,z_4\}$ be the corresponding configuration
in $\C \cup \infty$.  We rotate $S^2$ so that
\begin{equation}
\label{normalize}
z_4=\infty; \hskip 30 pt z_0 \in \R_+ \hskip 40 pt
E(z_0,z_4)=\max_{i<j} E(z_i,z_j).
\end{equation}

\noindent
{\bf The TBP:\/}
Now we discuss what the TBP looks like with this normalization.
The TBP has two kinds of points.   We say that the
{\it polar points\/} are the two antipodal points in the configuration.
We say that the three remaining points are {\it equitorial\/}.
When $z_4$ corresponds to a polar point, we have the following
configuration, which is unique up to the permutation of $z_1,z_2,z_3$:
\begin{equation}
\label{min1}
z_0=1; \hskip 30 pt z_1=\exp(-2 \pi i/3); \hskip 30 pt
z_2=0; \hskip 30 pt
z_3=\exp(2 \pi i/3); \hskip 30 pt
\end{equation}
When $z_4$ corresponds to an equitorial point, we have
\begin{equation}
\label{min2}
z_0=1; \hskip 30 pt z_1=-i \sqrt 3/2; \hskip 30 pt z_2=-1 \hskip 30 pt
z_3=i \sqrt 3/2;
\end{equation}
Later on, we will use a symmetry argument to avoid having to deal with the
second of these configurations.
\newline
\newline
{\bf The Space of Minimizers:\/}
We check $3$ facts using Equation \ref{metric1}.
\begin{enumerate}
\item If $|z|=1$ then $\|\Sigma^{-1}(z)-\Sigma^{-1}(\infty)\|=\sqrt 2$.
\item If $|z|\geq 2$ then $\|\Sigma^{-1}(z)-\Sigma^{-1}(\infty)\|<1$.
\item If $|x| \geq 4$ then $\|\Sigma^{-1}(z)-\Sigma^{-1}(\infty)\|<1/2$.
\end{enumerate}

Suppose now that $Z$ is a minimizer for some power law potential.
It follows from Lemma \ref{coulomb}
 and Item 3 that $|z_j|<3$ for
$j=1,2,3$.   If follows from Item 1 and
Lemma \ref{config} that $|z_0| \geq 1$.  We conclude that
$z_0 \in [1,4]$.    It follows 
Lemma \ref{coulomb2} and
Item 2 that $|z_j|\leq 2$ for $j=1,2,3$.   
We conclude that all the configurations we need to consider
lie in the compact region
compact region $\Omega=[0,4] \times [-2,2]^6 \subset \R^7$.
The map is given by just stringing out the coordinates of the points
$z_0,z_1,z_2,z_3$.
\newline
\newline
\noindent
{\bf Dyadic Objects:\/}
Given squares $Q_1, Q_2 \subset \C$, we write
$Q_1 \to Q_2$ if $Q_2$ is one of the $4$ squares
obtained by dividing $Q_1$ in half along both directions.
We say that a square $Q$ is a {\it dyadic square\/} if
there is a finite
\begin{equation}
[-2,2]^2 \to Q_1 \to \ldots \to Q_n=Q.
\end{equation}
The sides of $Q$ are necessarily parallel to the coordinate axes,
and the vertices have dyadic rational coordinates.  We call the
dyadic square {\it normal\/} if it does not cross the
coordinate axes.   A single subdivision of $[-2,2]$ produces
normal dyadic squares.

We make
the same definition for line segments as for squares., except
that the notation $S_1 \to S_2$ means that $S_2$ is one
of the two segments obtained by cutting $S_1$ in half.  We say
that a {\it dyadic segment\/} is a line segment $S$ such that
there is a finite chain
\begin{equation}
[0,4] \to S_1 \to \ldots \to S_n=S.
\end{equation}
We say that a {\it dyadic box\/} in $\Omega$ is a set of the form
$Q_0 \times Q_1 \times Q_2 \times Q_3$, where $Q_0$ is a dyadic
segment and $Q_j$ is a dyadic square for $j=1,2,3$.   The
whole space $\Omega$ itself is a dyadic box.
\newline
\newline
{\bf Subdivision:\/}
Let $B_j=Q_{j0} \times Q_{j1} \times Q_{j2} \times Q_{j3}$ be a
dyadic box for $j=1,2$.   We write
$B_1 \to_k B_2$ if
\begin{itemize}
\item $Q_{1i}=Q_{2i}$ for $i \not = k$.
\item $Q_{1k} \to Q_{2k}$.
\end{itemize}
The $k${\it th\/} {\it subdivision\/} of $B_1$ is the
union of all $B_2$ such that $B_1 \to_k B_2$.  When
$k=0$, this union consists of two dyadic boxes.  When
$k =1,2,3$, the union consists of $4$ dyadic boxes.
The set of all dyadic boxes in $\Omega$ forms a directed
tree.  Each dyadic box points to $14=2+4+4+4$ smaller
dyadic boxes.  

The divide-and-conquer algorithm will be
a depth-first search through this tree.  It seems that
the speed of the program depends a lot on how we
do the subdivisions, so we will explain this in
detail in the next section.

\subsection{The Divide and Conquer Algorithm}
\label{discuss1}

Our discussion applies to a general power law potential,
but we only apply the algorithm to the Coulomb potential.

Let $\epsilon>0$ be some small number.  In this section, we
explain in the abstract how we show, with a finite calculation,
that any winning configuration lies within $\epsilon$ of
one of the two TBP configurations described in \S \ref{space}.
There are $5$ components to our program:
\newline
\newline
{\bf Confinement:\/}
We say that a dyadic box 
\begin{equation}
{\cal Q\/}=Q_0 \times Q_1 \times Q_2 \times Q_3.
\end{equation}
 is $\epsilon$-{\it confined\/} if, for each $j=0,1,2,3$, the
dyadic square $B_j$ is contained in the open square of side length
$\epsilon$ centered on the point
$z_j$ from Equation \ref{min1}.
\newline
\newline
{\bf Tetrahedral Eliminator:\/}
We eliminate ${\cal Q\/}$ is if satisfies the hypotheses of
Lemma \ref{tetra eliminate}.  In this situation,
every configuration
in $\cal Q$ violates Equation \ref{tetra3}.
This function speeds up our calculations quite a bit.
\newline
\newline
{\bf Redundancy Eliminator:\/}
We will isolate a subset $\Omega' \subset \Omega$ with the
following property:  If $Z \in \Omega-\Omega'$,
then there is some $Z' \subset \Omega'$ (obtained
by permuting the points of applying some obviously
energy-decreasing move) such that $E(Z') \leq E(Z)$.
In \S \ref{redundancy} we will
describe some simple tests we use to show that
${\cal Q\/} \in \Omega-\Omega'$.   We eliminate
$\cal Q$ if it passes one of these tests.

One convenient property of $\Omega'$ is that
it contains the configuration in Equation \ref{min1}
but not the configuration in Equation \ref{min2}. This
means that we can automatically eliminate configurations
near the one in Equation \ref{min2}.  This is why our
notion of $\epsilon$-confinement only mentions
the configuration in Equation \ref{min1}.
\newline
\newline
{\bf Energy Estimator:\/}
Now we come to the main point.
Say that an {\it energy estimator\/} is a function
$\Phi: {\cal S\/} \to \R$ such with the following
property.
For every configuration
$Z=\{z_0,...,z_4\}$ with $z_j \in Q_j$, we have
$$E(Z) \geq \Phi({\cal Q\/}).$$  See Theorem
\ref{EE} for the definition of our Energy Estimator.
Theorem \ref{EE} is our main technical result.
\newline
\newline
{\bf Depth First Search:\/}
Recall that $M_e$ is the energy of the TBP.
Our program maintains a list $\cal L$ of dyadic boxes.
Initially, $\cal L$ just has the single box $\Omega$, the
whole space.   At a given stage of the program, the
algorithm examines the last box $\cal Q$ and eliminates it
if one of three things happens:
\begin{enumerate}
\item $\cal Q$ is $\epsilon$-confined.
\item Lemma \ref{tetra elim} eliminates $\cal Q$.
\item Onf the the tests in \S \ref{redundancy} eliminates $\cal Q$.
\item  $\Psi({\cal Q\/})>M_e$,.
\end{enumerate}
Otherwise, $\cal Q$ is eliminated and to $\cal L$ we append
the dyadic boxes in the $k$th subdivision of $\cal Q$ for
some $k \in \{0,1,2,3\}$.  We will explain how $k$ is
determined momentarily.
The algorithm halts if $\cal L$ is the empty list.  In this
case, we have shown that, up to symmetries,
 any minimizer lies within $\epsilon$ of the TBP.
\newline
\newline
{\bf Subdivision Rule:\/} The error term in Theorem \ref{EE} has the
form
\begin{equation}
\sum_{i=0}^3 \sum_{j\not = i} \epsilon(Q_i,Q_j).
\end{equation}
Here $\epsilon(Q_i,Q_j)$ is a quantity that depends on the geometry
of $Q_i$ and $Q_j$.   Roughly, it varies quadratically with the
side length of $Q_i$.   
We write
\begin{equation}
\epsilon(i)=\sum_{j \not = i} \epsilon(i,j).
\end{equation}
Then, again, $\epsilon(i)$ depends roughly quadratically on the
side length of $Q_i$.   We find the index $k$ that maximizes
the function $i \to Q(i)$ and then we use the $k$th subdivision
rule.    Thus, we subdivide in such a way as to try to make the
error estimate in Theorem \ref{EE} as small as possible.
When we compare this method with a more straightforward
method of subdividing so as to keep the squares all about
the same size, our method leads to a vastly faster computation.
\newline
\newline
\noindent
{\bf Remark:\/}
The main difficulty in our proof is choosing an energy
estimator that leads to a feasible calculation.  
We found the subtle Energy  Estimator from
Theorem \ref{EE} after quite a bit of trial and error.

\subsection{The Main Results}

Let ${\cal E\/}$ denote the energy function on $\Omega$.
For any $s>0$ let
$\Omega_s$ denote those configurations $\{z_k\}$
so that $z_k$ is contained in a square of side length $s$ centered
at the $k$th point of the TBP, normalized as in Equation \ref{min1}.
We shall be interested in the cases when 

\begin{equation}
s=2^{-11}.
\end{equation}

Running our computer program, we prove the following result.

\begin{lemma}[Confinement]
\label{maincomp}
Any Coulomb energy minimizer in $\Omega$ has the same
energy as some configuration in $\Omega_{s/4}$.
\end{lemma}

Let $H_e$ denote the Hessian
of $\cal E$, with respect to the potential
$E(r)=r^{-e}$.
In \S \ref{derproof} we will prove the following
result.

\begin{lemma}
\label{mainhess}
For $e = 1,2$ and
for any $W \in \Omega_{s}$, the matrix
$H_e(W)$ is positive definite.
\end{lemma}

Now let $Z_1$ be a $5$-point minimizer for the Coulomb potential.
By the Confinement Lemma, we can find a new configuration
$Z_2 \in \Omega_s$ with the same energy.   But $\cal E$ has
positive definite Hessian throughout $\Omega_s$, and
(by symmetry) the gradiant $\nabla {\cal E\/}$ vanishes at
$Z_0$, the TBP.   Restricting $\cal E$ to a straight line segment
connecting $Z_0$ to $Z_2$ we see that $\cal E$ is a
convex function with a local minimum at $Z_0$.  Hence
${\cal E\/}(Z_2)>{\cal E\/}(Z_0)$.  This proves
Theorem \ref{main}.
\newline
\newline
\noindent
{\bf Remarks:\/} 
\newline
(i) 
Our program also
establishes the Confinement Lemma for the function $E(r)=r^{-2}$. 
The same argument as above now establishes
Theorem \ref{main2}. \newline
(ii)
We didn't need to compute all the way down to $\Omega_{s/4}$.
We could have stopped at $\Omega_s$ in the Confinement Lemma.
However, an earlier version of this paper had a weaker result
on the Hessian, and
we did the extra computing to accomodate this.
There doesn't seem to be any reason to throw out our
stronger computaional result since we (or, rather, the computer)
took the trouble to get it.

\newpage

\section{Separation Estimates}

\subsection{Overview}
\label{quantities}

Let $\Sigma$ be stereographic projection.
The sets of interest to us have the form 
\begin{equation}
Q^*=\Sigma^{-1}(Q).
\end{equation}
where $Q$ is either a dyadic segment or a dyadic square
or the point $\infty$.  We will call $Q$ a
{\it dyadic planar set\/} and $Q^*$ a {\it dyadic spherical patch\/}.
A dyadic spherical patch is either the point $(0,0,1)$, or
an arc of a great circle, or else a subset of
the sphere bounded by $4$ arcs of circles.   

Let $(Q_1,Q_2)$ where  $Q_j$ is a dyadic planar set.
In this chapter we are interested in the following
two quantities.
\begin{equation}
\psi_{\rm max\/}(Q_1,Q_2)=\max \|p_1-p_2\|; \hskip 30pt p_j \in Q_j^*.
\end{equation}
\begin{equation}
\psi_{\rm min\/}(Q_1,Q_2)=\min \|p_1-p_2\|; \hskip 30 pt p_j \in {\rm Hull\/}(Q_j^*).
\end{equation}
Here ${\rm Hull\/}(Q_j^*)$ is the convex hull of $Q_j^*$.   For technical reasons
the lower bound we seek need to work for a slightly larger range of
points.  

In this chapter, we will estimate $\psi(Q_1,Q_2)$ in terms of quantities that can be
determined by a finite computation.   At the end of the chapter, we will
give, as an application of these estimates, the definition of the Tetrahedral
Eliminator discussed in \S \ref{discuss1}.

In the next definitions, we set $z=x+iy$.
If $Q$ is a dyadic square or dyadic segment, we define
$$
\underline x=\min_{z\in Q} |x|; \hskip 20 pt
\overline x=\max_{z\in Q} |x|; \hskip 20 pt
\underline y=\min_{z\in Q} |y|; \hskip 20 pt
\overline y=\max_{z\in Q} |y|; 
$$
\begin{equation}
\label{maxvertex}
\label{girth}
\delta=\frac{2 (\overline x-\underline x)}{1+\underline x^2+\underline y^2};
\hskip 30 pt \tau=\sqrt{{\rm dim\/}(Q)}.
\end{equation}
These quantities depend on $Q$ but we usually suppress
them from our notation.  
Here $\overline x-\underline x$ is just the sidelength of $Q$.
The quantity $\delta$ is an estimate on the side length of
$Q^*$.     Note that
$\tau=1$ if $Q$ is a dyadic segment and
$\tau=\sqrt 2$ if $Q$ is a dyadic square.

When $Q=\{\infty\}$ we
set $\delta=\tau=0$ and we don't define the other quantities
know $Q^*$ exactly in this case.

Given a pair $(Q_1,Q_2)$ of dyadic objects, we 
\begin{equation}
D=\|\Sigma^{-1}(z_1)-\Sigma^{-1}(z_2)\|; \hskip 30 pt
\delta=\frac{\delta_1\tau_1+\delta_2\tau_2}{4}.
\end{equation}
Here $z_j$ is the center of $Q_j$ and
$\delta_j=\delta(Q_j)$, as defined above.
Here is our main result.

\begin{lemma}[Bound]
\label{fine}
\label{bound}
The following is true relative to any pair of dyadic patches.
$$
\psi_{\rm min\/} \geq D(1-\delta^2/2)-(\sqrt{4-D^2}) \delta;
\hskip 20 pt
\psi_{\rm max\/} \leq D + (\sqrt{4-D^2}) \delta.$$
\end{lemma}

\noindent
{\bf Remark:\/} Our proof will yield the better bounds
$$\psi_{\rm min\/} \geq D \cos(\delta)-(\sqrt{4-D^2}) \sin(\delta)
\hskip 20 pt
\psi_{\rm max\/} \leq D \cos(\delta) + (\sqrt{4-D^2}) \sin(\delta),$$
but for  computational reasons we want to avoid the trig functions.
We have used rational replacements which are quite close in
practice to the trig functions.
\newline

There is one situation where we can get a completely sharp upper bound.
We define
\begin{equation}
\psi'_{\rm max\/}(Q_1,Q_2)=\max \|p_1-p_2\|; \hskip 30pt p_j \in Q_j^*.
\end{equation}
This time we mazimize over pairs $(p_1,p_2)$, where
$p_j$ is a vertex of $Q_j^*$.   A finite computation
gives $\psi'_{\rm max\/}$.  We define
$\psi'_{\rm min\/}$ similarly.
Recall that a dyadic square is
normal if it doesn't cross the coordinate axes.

\begin{lemma}[Perfect Bound]
\label{perfect}
If $Q_1$ is normal and
$Q_2=\{\infty\}$ then we have
$\psi_{\rm max\/}=\psi'_{\rm max\/}$
and
$\psi_{\rm min\/}=\psi'_{\rm min\/}.$
\end{lemma}

We define $\Psi_{\rm min\/}$ to be the best bound we can
get from the two lemmas above (or $0$, if no lemma applies.)
We define $\Psi_{\rm max\/}$ to be the best bound we can
get from the two lemmas above (or $2$, if no lemma applies.)
\newline

The rest of the chapter is devoted to proving Lemma
\ref{fine} and Lemma \ref{perfect}.

\subsection{Proof of Lemma \ref{fine}}

\begin{lemma}
\label{arcs}
Let $A$ and $B$ be arcs of the unit circle.
Let $\Delta_A$ and
$\Delta_B$ denote the arc lengths of $A$ and $B$ 
Let $D_A$ and $D_B$ denote
the distance between the endpoinds of $A$ and $B$ respectively.
Let $$\delta=\frac{\Delta_A-\Delta_B}{2}.$$
Then
$$D_B = D_A \cos(\delta) \pm \bigg(\sqrt{4-D_A^2}\bigg)\sin(\delta).$$
\end{lemma}

\startproof
We have the relations
\begin{equation}
D_A=2 \sin(\Delta_A/2); \hskip 30 pt
D_B=2 \sin(\Delta_B/2)=2\sin(\Delta_A/2-\delta).
\end{equation}
Hence, by the angle addition formula,
\begin{equation}
\label{anglesum}
D_B=D_A \cos(\delta) - 2\cos(\sin^{-1}(D_A/2))\sin(\delta).
\end{equation}
Noting that
\begin{equation}
\cos(\sin^{-1}(x))=\sqrt{1-x^2},
\end{equation}
 we see that Equation
\ref{anglesum} is the same as the result we want.
\endproof

\begin{lemma}
\label{confine}
Let $Q$ be a dyadic set.
Every point of $Q^* $ is within $\delta \tau/2$ units of $\Sigma^{-1}(z_Q)$.
\end{lemma}

\startproof
When evaluated at all points of $Q$, the quantity in
Equation \ref{metric2} is maximized at the point $(\underline x,\underline y)$.   Its
value at that point is exacrly 
$$\frac{2}{1+\underline x^2+\underline y^2} = \frac{\delta}{s}.$$
Therefore $\Sigma^{-1}$ expands distances on $Q$ by at most a factor
of $\delta/s$, and every point of $Q$ is within $s\tau/2$ of $z$.
From here the result is obvious.
\endproof

Now we give the main argument for Lemma \ref{fine}.  
We consider the lower bound first.
Given a point $p \in S^2$ and some $r>0$ we let
$H(p,\epsilon) \subset S^2$ denote the convex hull
of the set of points on $S^2$ that are within
$\epsilon$ of $p$ in terms of arc length on $S^2$.
We call such a set a {\it cap\/}.

It follows from Lemma \ref{confine} and the convexity
of caps that
\begin{equation}
{\rm Hull\/}(Q_j^*) \subset H_j=H(p_j,\delta_j\tau_j/2).
\end{equation}
Here $p_j=\Sigma^{-1}(z_j)$, where $z_j$ is the
center of $Q_j$.   

Suppose first that $H_1$ and $H_2$ are disjoint.
Let $(q_1,q_2) \in H_1 \times H_2$ be
two points which realize the minimum of
$\|q_1-q_2\|$.    We must have $q_j \in \partial H_j$.
Also, the segment joining $q_1$ to $q_2$ must be
perpendicular to both $\partial H_1$ and $\partial H_2$.
This situation leads to the result that $q_1$ and $q_2$ are
contained on the great circle $C$ joining $p_1$ and $p_2$.

We can now reduce everything to a problem in the plane.
Let $\Pi$ be the plane containing $C$.   We identify
$\Pi$ with $\C$, so that $C$ is the unit circle.
Figure 3.1 shows the situation.  We have drawn the case
when both intersections lie in a half-disk, but this feature
is not a necessary part of the proof.

\begin{center}
\resizebox{!}{2.6in}{\includegraphics{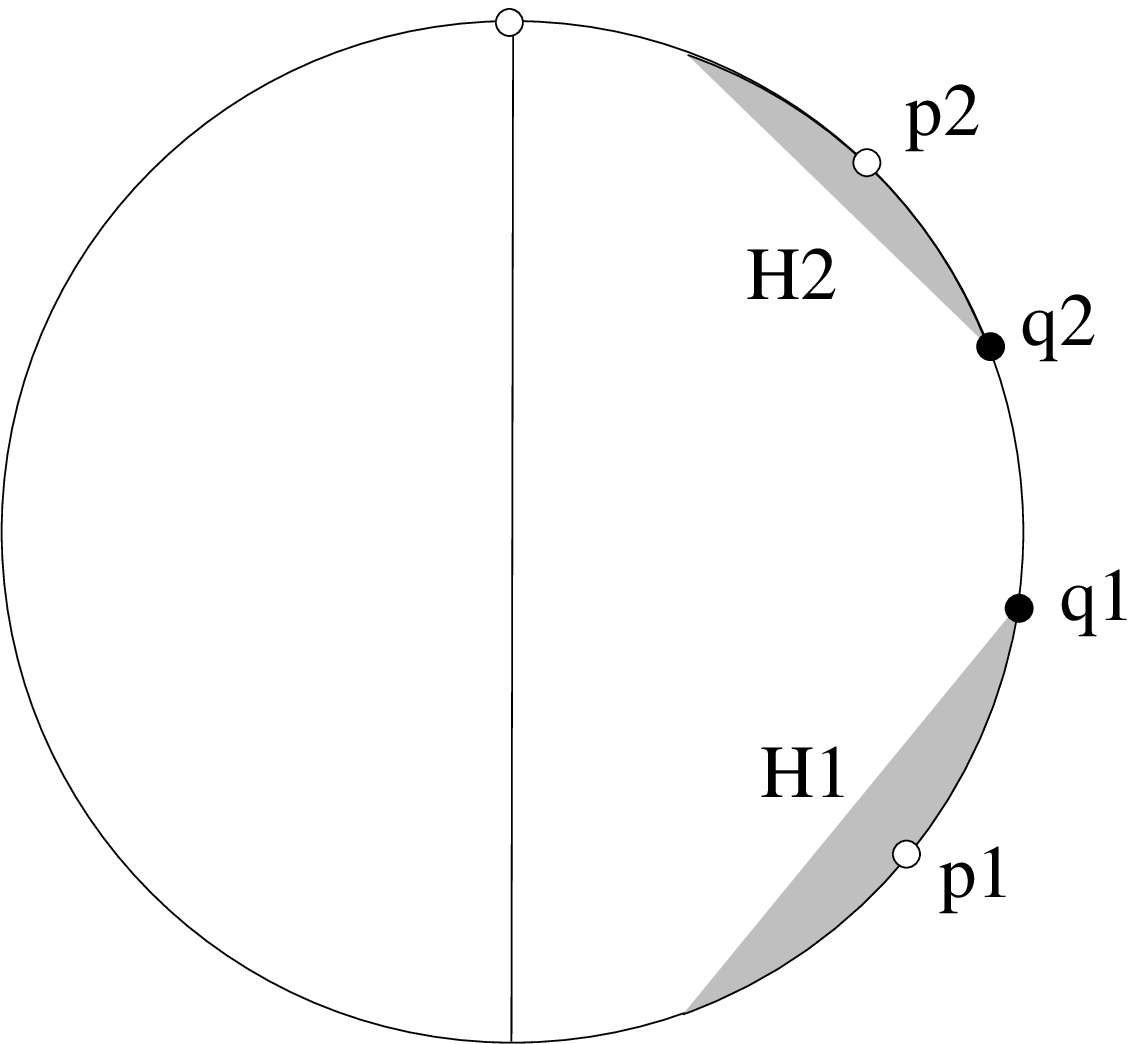}}
\newline
{\bf Figure 3.1:\/} Two circular caps.
\end{center}

Let $A$ be the short circular arc joining $p_1$ and $p_2$, and let
$B$ be the short circular arc joining $q_1$ and $q_2$.
Referring to Lemma \ref{arcs}, we have
\begin{equation}
D_A=D(Q_1,Q_2)=\|p_1-p_2\|;
\end{equation}
Here $D=D(Q_1,Q_2)$ is the quantity in the statement of the lemma.
We also have
\begin{equation}
\psi_{\rm min\/} \leq D_B=\|q_1-q_2\|.
\end{equation}
Finally, the length of the circular arc joining $p_j$ to $q_j$ is
$\delta_j\tau_j/2$.  Therefore
\begin{equation}
\delta(A,B)=\Delta_A-\Delta_B=
\frac{\delta_1\tau_1/2+\delta_2\tau_2/2}{2}=
\delta(Q_1,Q_2).
\end{equation}
Applying Lemma \ref{arcs}, and using our three equations, we get
\begin{equation}
\psi_{\rm min} \geq D \cos(\delta) - \bigg(\sqrt{4-D^2}\bigg) \sin(\delta)
\end{equation}
But $\cos(\delta) \geq 1-\delta^2$ and
$\sin(\delta) \leq  \delta$.
This gives us the lower bound from Lemma \ref{fine} in the case
the caps are disjoint.

When the caps intersect, we have $\psi_{\rm min\/}=0$.  We just want to show that
the lower bound in Lemma \ref{fine} is nonpositive.     We want to reduce this
case to the previous one.   Choose some small $\epsilon>0$ and consider smaller
caps $H_1'$ and $H_2'$ that are separated by a distance of exactly
$\epsilon$.   These two caps are based on some number $\delta'<\delta$.
Our argument in the previous case gives
\begin{equation}
\epsilon>D(1-(\delta')^2/2)-\bigg(\sqrt{4-D^2}\bigg) \delta'
\end{equation}
\begin{equation}
D(1-(\delta')^2/2)-\bigg(\sqrt{4-D^2}\bigg) \delta'>
D(1-\delta^2/2)-\bigg(\sqrt{4-D^2}\bigg) \delta.
\end{equation}
Combining these two results, we see that the lower bound in
Lemma \ref{fine} is less than $\epsilon$.  But $\epsilon$ is
arbitrary. Hence, the lower bound in Lemma \ref{fine}
is nonpositive in this case.
This completes the proof of the lower bound.
\newline

The proof of the upper bound is similar. Suppose first that $H_1 \cap H_2$ does
not contain a pair of antipodal points.  Then, by Lagrange Multipliers,
the pair of points $(q_1,q_2)$ realizing the maximum lie on the great
circle through $p_1$ and $q_2$.   We then rotate as above and
apply the same argument.  This time we get
\begin{equation}
\psi_{\rm max} \leq D \cos(\delta) + \bigg(\sqrt{4-D^2}\bigg) \sin(\delta)
\end{equation}
But $\cos(\delta) \leq 1$ and $\sin(\delta) \leq \delta$.
This gives us the lower bound from Lemma \ref{fine} in the case
$H_1 \cup H_2$ does not contain a pair of antipodal points.

If $H_1 \cup H_2$ contains a pair of antipodal points then one of two things is
true.  If $p_1$ and $p_2$ are antipodal points, then obviously the upper bound
in Lemma \ref{fine} gives a number greater than $2$.   If $p_1$ and $p_2$ are
not antipodal we can use the same shrinking trick that we used in the
previous case to show that the upper bound in Lemma \ref{fine} gives
a number that is at least $2$.  In either case, the upper bound still holds.

This completes the proof of Lemma \ref{fine}.

\subsection{Proof of Lemma \ref{perfect}}
\label{perfectproof}

\begin{lemma}
\label{perfect2}
Suppose that $Q_1$ is a dyadic segment or square and
$Q_2=\{\infty\}$.  Let $p_2=(0,0,1)$.  Then, for point
$p_1 \in Q_1^*$ we have
$\psi'_{\rm min\/} \leq \|p_1-p_2\| \leq \psi'_{\rm max}.$
In particular, the upper  bound in Lemma \ref{perfect} is true.
\end{lemma}

\startproof
Consider the lower bound first.
 The point $z_1 \in Q_1$ which
minimizes 
\begin{equation}
\label{disks}
\|\Sigma^{-1}(z_1)-\Sigma^{-1}(\infty)\|
\end{equation}
 is the
point furthest from the origin.  But, since disks are convex, the
point of $Q_1$ farthest from the origin is a vertex.

Now for the upper bound.  The point of $Q_1$ that maximizes
Equation \ref{disks} is the one closest to the origin.
If a vertex of $Q_1$ does not minimize the distance to
the origin, then (by Lagrange multipliers) some ray through
the origin intersects a side of $Q_1$ at a right angle.  Since
the sides of $Q_1$ are parallel to the coordinate axes,
this can only happen if $Q_1$ crosses one of the coordinate
axes.  But $Q_1$ does not
cross the coordinate axes.
\endproof

Lemma \ref{perfect2} immediately gives the upper bound.
For the lower bound (which makes a different statement)
we need to deal with all the points in the convex hull
$H_1={\rm hull\/}(Q_1^*)$.
Any point $q \in H_1$ that minimizes
$\|q-(0,0,1)\|$ must lie on $\partial H_1$.  We claim that
$\partial H_1$ is the union of the following
\begin{itemize}
\item $Q_1^*$.
\item The flat quadrilateral $F$ that is the convex hull of the vertices of $Q_1^*$.
\item The convex hulls $H(E_j)$ of the edges $E_1,...,E_4$ of $Q_1^*$.
\end{itemize}
Each of these sets is a subset of $H_1$ and we easily check, in each
case, that every point in each set is a boundary point.   Since the
union of these sets is a topological sphere, it must account for the
entire boundary.

Our lemma above takes care of the case when $q \in Q_1^*$.   If
$q$ is an interior point of $F$, then the segment joining $(0,0,1)$
to $q$ is perpendicular to $F$.  But $F$ is completely contained in a
hemisphere that has $(0,0,1)$ on its boundary, so no such
segment can exist.  

Finally, suppose that $q \in H(E_j)$.  Note that $H(E_j)$ is contained
in a plane $\Pi$ which also contains $(0,0,1)$.  But, considering
the distance minimizimation problem in $\Pi$, we see that
$q$ cannot be an interior point of $H(E_j)$.  (The picture
looks just like the one drawn in Figure 3.1.)   But then 
either $q \in Q_1^*$ or $q \in F$, the two cases we have
already handled.    This completes the proof.

\subsection{The Tetrahedral Eliminator}
\label{tetra4}

Let $\cal Q$ be a dyadic box, as in the previous section.
As in Equation \ref{tetra3}, the quantity $M_E$ is the
energy of the TBP and the quantity $T_E$ is the energy of
the regular tetrahedron.   Let 
$\Psi_{\rm max\/}$ be the upper bound on the distance
between $Q_1^*$ and $Q_1^*$ as in \S \ref{quantities}.

\begin{lemma}[Tetrahedral Eliminator]
\label{tetra eliminate}
\label{tetra elim}
Let $\cal Q$ be a dyadic box.  
Suppose that there is some index $i$ such that
\begin{equation}
\label{tetra5}
\sum_{j \not=i} E(d_{ij}) > M_E-T_E;
\hskip 30 pt d_{ij}=\Psi_{\rm max\/}(Q_i,Q_j).
\end{equation}
Then no configuration in $\cal Q$ is a minimizer for the $E$-energy.
\end{lemma}

\startproof
By construction 
$\|p_i-p_j\| \leq d_{ij}$
 for all
$p_i \in Q_i^*$ and $p_j \in Q_j^*$.
Therefore
\begin{equation}
E(P,p_i)>\sum_{j \not = i} E(d_{ij})>M_E-T_E
\end{equation}
for any configuration $P$ corresponding to a point in $\cal Q$.
But then every such configuration violates Equation \ref{tetra3}
and cannot be a minimizer.
\endproof

\subsection{Discussion}

We can get a cheap Energy Estimator using the fact that
\begin{equation}
{\cal E\/}(Z) \geq  \sum_{i<j} E(d_{ij}); \hskip 30 pt
d_{ij}=\Psi_{\rm max\/}(Q_i,Q_j).
\end{equation}
However, experiments lead us to
believe that the main calculation would 
effectively take forever using this estimator.

Here is the problem:  Even though $E(d_{ij})$ gives
a good estimate for the minimum energy of a pair
of points $(p_i,p_j) \in Q_i^* \times Q_j^*$, there
is no guarantee that we can find a \underline{single}
configuration $\{p_j\}$ such that each $d_{ij}$ is
nearly realized by $E(\|p_i-p_j\|)$.   So, 
the minimum energy of a configuration we can
actually produce might be much higher than
our estimate.  That is, our estimate might not
be that good {\it globally\/} (for all $10$ interactions)
even though it is good {\it locally\/} (for pairs of
interactions.)  
We need to work harder to get a globally good energy
estimator.  This is the content of Theorem \ref{EE}.

\newpage

\section{The Redundancy Eliminator}
\label{redundant0}

The constructions in this chapter work for any power law,
and most of the constructions (just symmetry) work for
any decreasing energy function.

\subsection{Inversion}

An {\it inversion\/} is a conformal involution $\rho$  of $\C \cup \infty$ that
fixes some circle $C$ pointwise. We call $\rho$ {\it stereo-isometric\/} if $\Sigma^{-1}(C)$ is
a great circle of $S^2$.   This condition means that
stereographic projection conjugates $\rho$ to an isometry
of $S^2$.
The purpose of this section is to prove the following lemma.
The reader might want to just note the result and skip
the proof on the first reading.

\begin{lemma}
\label{invert}
Let $\rho$ be a stereo-isometric inversion fixing
a circle centered on a point of $[1,\infty) \subset \R$.
Let $X \subset \C$ be a strip bounded by horizontal
lines, one lying above $\R$ and one lying below $\R$.
Let $X_- \subset X$ denote those points $z=x+iy$ with
$x<1-\sqrt 2$.    Then $\rho(X_-) \subset {\rm interior\/}(X)-X_-$.
\end{lemma}

\startproof
Being an inversion, the map $\rho$ has $3$ nice geometric properties.
\begin{itemize}
\item $\rho$ interchanges the disk $D$ bounded by $C$ with its complement.
\item $\rho$ interchanges the center of $C$ with $\infty$.
\item $\rho$ preserves each ray through the center of $C$.
\end{itemize}
We will use these properties in our proof.

Let $z_0$ be the center of $C$.  We claim that $C$
separates $z_0$ from any point on
$\R$ that lies to the left of $1-\sqrt 2$.
Let $p_0=\Sigma^{-1}(z_0)$.  As usual,
$(0,0,1)=\Sigma^{-1}(\infty)$.   Let $\Pi$
be the plane equidistant between $p_0$ and $(0,0,1)$.
Then 
\begin{equation}
C=\Sigma(\Pi \cap S^2)
\end{equation}
Let $H$ be the half-space bounded by $\Pi$ that contains $p_0$.

The point $p_0$ lies on the great circle
connecting $\Sigma^{-1}(1)=(1,0,0)$ to $(0,0,1)$.   Since $z_0 \in [1,\infty)$,
the point $p_0$ lies on the short arc connecting $(1,0,0)$ to $(0,0,1)$.
But then 
\begin{equation}
q=\Sigma^{-1}(1-\sqrt 2) =
\bigg(-\frac{1}{\sqrt 2},0,-\frac{1}{\sqrt 2}\bigg) \not \in {\rm interior\/}(H).
\end{equation}
The extreme cases occurs when $p_0=(1,0,0)$. In this extreme case
$q \in \Pi \cap S^2$.  As $p_0$ moves towards $(0,0,1)$,
$H \cap S^2$ moves away from $q$.   Since $H$ does not
contain $(0,0,1)$ either, the entire arc connecting $q$ to
$(0,0,1)$ is disjoint from the interior of $H$. But this means
that $C$ separates $z_0$ from all points of $\R$ that lie
to the left of $1-\sqrt 2$.
This proves our claim.

Let $D$ be the
disk bounded by $C$ that contains $z_0$.
From what we have already shown,
the interior of $D$ is disjoint from the vertical line $L$
through $1-\sqrt 2$.  Therefore
$D \cap X_-=\emptyset$.
Since $\rho$ swaps $D$ and its complement, we have

\begin{equation}
\rho(X_-) \subset D.
\end{equation}
  Hence
\begin{equation}
\label{contain1}
\rho(X_-) \cap X_-=\emptyset.
\end{equation}

At the same time, $\rho$ preserves all rays through the $z_0$,
the center of $D$.   Noting again that $z_0 \in \R$, we see
geometrically that $\rho(w)$ is closer to $\R$ than $w$
for all $w \in X-D$.   See Figure 4.1.
In particular, we have
\begin{equation}
\label{contain2}
\rho(X_-) \subset {\rm interior\/}(X).
\end{equation}
The first statement of the lemma now follows from Equations \ref{contain1}
and \ref{contain2}.  The second statement follows from the first
statement and from the fact that $\rho$ is an involution.
\endproof

\begin{center}
\resizebox{!}{2.3in}{\includegraphics{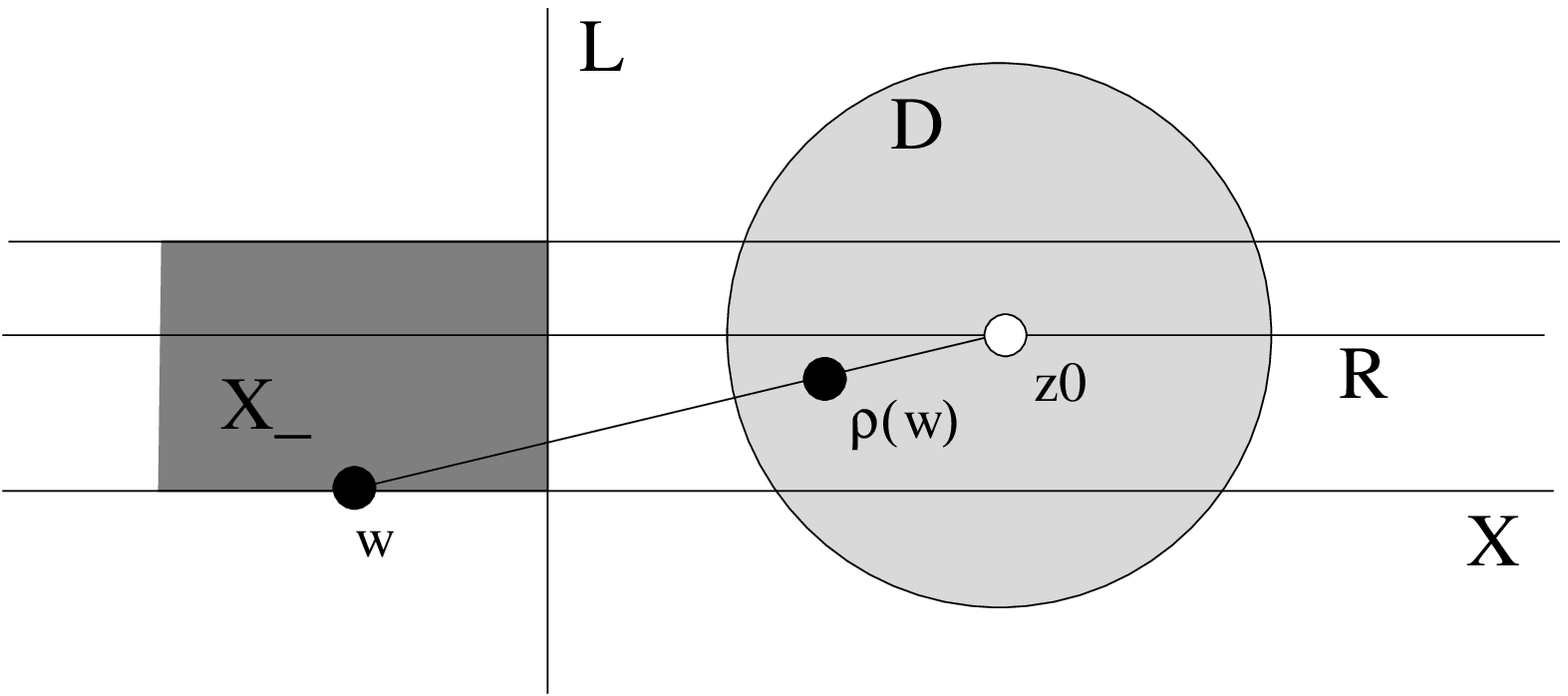}}
\newline
{\bf Figure 4.1:\/} Inversion
\end{center}

\subsection{A Lemma about Five Points}

The result in this section is certainly well known.
We suggest that perhaps the reader just
note the result and skip the proof on the first reading.

\begin{lemma}
\label{mindist2}
\label{config}
Let $P$ be any configuration of $5$ distinct points on $S^2$.  Then
some pair  points of $P$ are within $\sqrt 2$ units of each other.
\end{lemma}

\startproof
A {\it quarter sphere\/} $B$ is a
region on $S^2$ bounded by two semicircular great arcs,
meeting at antipodal points, such that the interior angles at
the intersection points are $\pi/2$.   
The {\it axis\/} of $B$ is the circular arc,
contained in $B$, that connects the midpoints of
the two bounding arcs.   We will use the following
easy-to-prove principle:  If $p$ is any
point on the axis of $B$, then the closed hemisphere
centered at $p$ contains $B$.   Call this
{\it Property X\/}.

The lemma has a trivial proof if $P$ contains two antipodal
points, so we assume this is not the case.
Let $H_j$ be the closed
hemisphere centered at $p_j$ of $P$.
Another way to state this lemma is that there are indices
$i \not = j$ such that $p_i \in H_j$.
Assume this is false, for the sake of contradiction.

Let $N$ denote the hemisphere centered at $(0,0,1)$, and
let $S$ denote the opposite hemisphere.
We normalize so that $p_0=(0,0,1)$.  By
assumption, $p_1,...,p_4 \in S-N$.   Since $P$
does not contain antipodal points, $P$ does not contain
$(0,0,-1)$.  But then,
for $j=1,2,3,4$, 
there is a unique point $p_j^*$ on the equator $N \cap S$
which is is as close as possible to $p_j$.

Let $H_j^*$ denote the hemisphere centered 
centered at $p_j^*$.   By construction, $H_j^* \cap S$ is
a quarter sphere and $p_j$ is a point on its axis.  By Property X,
we have $H_j^* \cap S \subset H_j$.      But then,
by assumption, 
$p_i \not \in H_j^* $ if $i \not = j$.
Let $\pi: S \to \R^2$ be the projection map
$\pi(x,y,z)=(x,y)$.  Then $\pi(S)$ is the unit disk
and $\pi(S \cap H_j^*)$ is a half disk.
Finally, $\pi(p_i)$ and $\pi(p_i^*)$ lie on the
same ray through the origin.  For this reason,
$p_i \in H_j^*$ if and only if $p_i^* \in H_j^*$.
Since $p_i \not \in H_j^*$, we conclude that
$p_i^* \not \in H_j^*$.

Now we know that $p_i^* \not \in H_j^* $ if $i \not = j$.
Hence, the distance from $p_i^*$ to $p_j^*$
is greater than $\pi/2$.   But then we have
$4$ points $p_1^*,...,p_4^*$, contained on the
equator, each of which is more than $\pi/2$ from
any of the other points.  This is a contradiction.
\endproof

\subsection{The Set of Good Configurations}
\label{good}

Recall that $\Omega$ is our configuration space. 
Our construction here works for any energy function $E$.
In this section we define the set $\Omega' \subset \Omega$, as discussed
in \S \ref{discuss1} in connection with the Redundancy Eliminator.

Let $Z=\{z_0,...,z_4\}$ be a configuration of points in $\C \cup \infty$,
normalized as in \S \ref{twomin}.   We write
$z_j=x_i+y_i$.  Also, we set $p_j=\Sigma^{-1}(z_j)$.
Let $\Omega'$ be the set of configurations $Z$ having all the
following properties.

\begin{enumerate}
\item $\|p_i-p_j\| \leq \|p_0-p_4\|$ for all indices $i \not = j$.  Also $x_0 \geq 1$.
\item $y_1 \leq 0 \leq y_2 \leq y_3$.
\item If $y_1<0<y_3$ then $x_2 \geq 1-\sqrt 2$.
\end{enumerate}

\begin{lemma}
\label{redundant}
For any $Z \in \Omega$, there exists 
$Z' \in \Omega'$ such that
$E(Z') \leq E(Z)$.
\end{lemma}

\startproof
Permuting the points, we can the first half of Property 1.
By Lemma \ref{config} and the first half of
Property 1, we must
have $\|p_0,p_4\| \leq \sqrt 2$.    But this
fails if $z_0 \in [0,1)$.  Hence $x_0 \geq 1$.
In short, the first half of Property 1 implies the
second half.

Reflecting in $\R$, we can arrange that $y_2 \geq 0$.
Permuting again, we can retain Property 1 and arrange
that $y_1 \leq y_2 \leq y_3$.  To get Property 2, we just
have to deal with the situation when $0<y_1$.
Let $Z'$ be the new configuration obtained when we
replace $z_1$ with the conjugate $\overline z_1$ and otherwise
keep the points the same.   Let $P$ and $P'$ be the
corresponding configurations in $S^2$.   The
points $p_0,...,p_5$ are all contained in the hemisphere
bounded by the great circle 
\begin{equation}
C=\Sigma^{-1}(\R \cup \infty).
\end{equation}
The configuration $P'$ is obtained from $P$ simply by
reflecting $p_1$ across $C$.  But, as can be seen from e.g.
the Pythagorean theorem, we then have
$\|p_1'-p_j\| \geq \|p_1-p_j\|$ for $j \not 3$.
The other distances do not change.  Since $E$ is
a monotone decreasing function, we have
$E(Z') \leq E(Z)$.  This gives us Property 2.

It remains to deal with Property 3.   
Suppose that $y_1<0<y_3$ and $x_2<1-\sqrt 2$.
   Let $\rho$ be the stereo-isometric
inversion that swaps $z_0$ and $z_4$.
Let $Z'=\rho(Z)$.    Let $z_j'=\rho(z_j)$.
We set $z_j'=x_j'+iy_j'$.   By construction
$Z'$ has property 1.
By symmetry, 
\begin{equation}
\label{si1}
y_1'<0<y_3'; \hskip 30 pt y_2' \geq 0.
\end{equation}
There are two cases to consider.
\newline
\newline
{\bf Case 1:\/}
Suppose first that
\begin{equation}
\label{si2}
y_2' \leq y_3' 
\end{equation}
Equations \ref{si1} and \ref{si2} combine to say that
$Z'$ has Property 2 as well.
Let $X$ be the strip bounded by the horizontal lines through $z_1$ and $z_3$.
Then $X$ satisfies the hypotheses of Lemma \ref{invert}.  Moreover,
$z_2 \in X_-$, the half-strip from Lemma \ref{invert}.
Then $z_2' \in X-X_-$ by Lemma \ref{invert}. Hence $x_2'>1-\sqrt 2$.
This shows that $Z'$ has Property 3.  Hence $Z' \in \Omega'$.
\newline
\newline
{\bf Case 2:\/}
Suppose that $y_2'>y_3'$.  Then we let $Z''$ be the configuration
obtained from $Z'$ by swapping $z_2'$ and $z_3'$.   Then
$Z''$ satisfies Properties 1 and 2.  The argument in
Case 1 again shows that
$z_2' \in X-X_-$.
In particular, $z_2' \in X$.
By construction, 
\begin{equation}
z_3 \not \in {\rm interior\/}(X).
\end{equation}
By Lemma \ref{invert}, we have
\begin{equation}
\rho(X_-) \subset {\rm interior\/}(X).
\end{equation}
Combining these last two equations, we see that
\begin{equation}
z_3 \not \in \rho(X).
\end{equation}
Since $\rho$ is an involution, this last equation gives
us
\begin{equation}
z_3' \not \in X_-.
\end{equation}
Summarizing the situation, we now know that
\begin{equation}
z_2' \in X; \hskip 30pt z_3' \not \in X_-.
\end{equation}

Since $0<y_3'<y_2'$ and $z_2' \in X$, we have
\begin{equation}
z_3' \in X.
\end{equation}
Combining the last two equations, we have
\begin{equation}
z_3' \in X-X'; \hskip 30 pt \Longrightarrow \hskip 30 pt
x_2''=x_3' \geq 1-\sqrt 2.
\end{equation}
The first statement implies the second.  
The second statement shows that $Z''$ has Property 3 as well.
Finally, $E(Z'')=E(Z)$.
\endproof

\subsection{The Main Construction}
\label{redundancy}

Now we define the Redundancy Eliminator discussed in \S \ref{discuss1}.
The redundancy eliminator performs $4$ tests, one per property
discussed above.  
Let ${\cal Q\/}=Q_0 \times Q_1 \times Q_2 \times Q_3$ be a dyadic box.
As usual, a configuration $Z \in \cal Q$ defines points
$p_0,...,p_4$, with $p_j \in Q_j^*$.
\newline
\newline
{\bf Property 1:\/}
Let $\Psi_{\rm max\/}$ and
$\Psi_{\rm min\/}$ be the separation functions from
\S \ref{quantities}.
We eliminate $\cal Q$ if
\begin{equation}
\Psi_{\rm max\/}(Q_i,Q_j)<\Psi_{\rm min\/}(Q_0,Q_4)
\end{equation}
for some pair of indices $i<j$ such that $(i,j) \not = (0,4)$.
We also eliminate $\cal Q$ if $\underline x_0<1$.
\newline
\newline
{\bf Property 2:\/}
Let $\underline x_j$, etc be as in \S \ref{quantities}.
We eliminate $\cal Q$ if any of the following is true.
\begin{itemize}
\item $\underline y_1 \geq \overline y_2$;
\item $\underline y_2 \geq \overline y_3$.
\item $\overline y_2 \leq 0$.
\item $\underline y_1\geq 0$.  
\end{itemize}
In the first $3$ cases, the fact that we are using a
weak inequality rather than a strict inequality means
that sometimes we eliminate some configurations
that lie in $\partial \Omega'$.   Since we want to
consider every configuration in $\Omega'$,
we need to justify this.  We give the justification in the
next section.
\newline
\newline
{\bf Property 3:\/} 
We eliminate $\cal Q$ if all of the following happen.
\begin{itemize}
\item $\overline y_1<0$.
\item $\underline y_3 >0$.
\item $\overline x_2<1-\sqrt 2$.
\end{itemize}
This time there are no boundary cases to worry about.

\subsection{Discussion of Boundary Cases}
\label{boundary}

Here we discuss the use of weak inequalities in connection
with Property 2 in the previous section.
First of all, the reason why we want to do this is that
it speeds up the computation.  For example, were
we to keep dyadic boxes with $\overline y_2=0$
we would need to consider many more dyadic
boxes that are near the TBP.   

The justification for why we can use weak
inequalities is that all the configurations in
$\partial \Omega'$ that we eliminate are
actually counted twice, and we do not eliminate
the relevant dyadic boxes both times.

To clarify the situation, we make the
interpretation that our dyadic
squares $Q_2$ and $Q_3$ are 
missing their top boundaries and $Q_1$ is missing
its bottom boundary. 
With this convention, the divide and
conquer algorithm from \S \ref{discuss1}
examines every point in the configuration
space except those for which $|y_j|=2$ for some $j$.
These configurations are not minimizers.
See the discussion in \S \ref{twomin}.
  The main point here is that the
union of ``quarter-open'' dyadic squares in the
subdivision of a ``quarter open'' dyadic square 
is still equal to the original ``quarter open''
dyadic square.  In the case of $Q_2$ and $Q_3$, the
bottom edges fill in for the top edges. In the case of
$Q_1$ the top edges fill in for the bottom ones.

With this interpretation, we can simply eliminate the
dyadic box $\cal Q$ mentioned above, because
it contains no configurations in $\Omega'$.
Adopting this convention has no effect whatsoever on
our program.  It is simply a question of how we interpret
the output of the program.

\newpage

\section{The Energy Estimator}
\label{energy}

\subsection{Preliminaries}

We fix some value of $n$.
Say that a {\it block\/} $\cal Q$ is a collection
$Q_0,...,Q_{n}$ of dyadic objects, with
$Q_{n}=\{\infty\}$ and all the other objects
either segments or squares.
Say that a collection of points $z_0,...,z_{n}$ of points
is {\it dominated by\/} the block if
$z_k \in Q_k$ for all $k$.
In our application, we will take $n=4$.
In this case, the set of configurations dominated by
a block is precisely a dyadic box.
\newline
\newline
{\bf Remark:\/}
Though we always take $Q_k$ to be a dyadic segment
or square for $k=0,...,n-1$, the
interested reader will note that everything we do
works the same way with the milder constraint that
each $Q_k$ is either a segment in $[0,\infty)$ or
a square with sides parallel to the coordinate axes
that does not cross the coordinate axes.
\newline

We say that $\{z_k\}$ is a {\it vertex configuration\/} if
$z_k$ is a vertex of $Q_k$ for $k=0,...,n$.    
There is a finite list of vertex configurations.  
Given an energy function $E$, we define
\begin{equation}
\label{maineng}
{\cal E\/}({\bf \cal Q\/})=\min_{Z} E(Z).
\hskip 30 pt
{\cal E\/}'({\bf \cal Q\/})=\min_{Z'} E(Z').
\end{equation}
The first minimum is taken over all configurations $Z$
dominated by $\cal Q$ and the second
minimum is taken over all vertex configurations $Z'$.
The purpose of our main result in this chapter is to bound the quantity
\begin{equation}
{\bf ERR\/}(\cal Q)=
{\cal E\/}'({\bf \cal Q\/})-
{\cal E\/}({\bf \cal Q\/})
\end{equation}
in terms of $\cal Q$.

We are interested mainly in the case when the function
$E$ is a power law, but we will state things more generally
just to clarify the formulas we get.
We suppose that $E$ is convex decreasing and
satisfies the following property.
 If $0<r_2 \leq r_1$ then there are
constants $c \geq 1$ and $h \geq 0$ (depending on everything)
 such that
\begin{equation}
\label{nice}
\frac{E'(r_2)}{r_2}=c \frac{E'(r_1)}{r_1}; \hskip 30 pt
E''(r_2)=(c+h) E''(r_1).
\end{equation}
This condition could be stated more simply, but we have
stated exactly the version we will use in Lemma \ref{comparison1},
the one place where we use it.  Other than this
one place, the rest of our argument works for any
convex decreasing function.
The power law functions all
satisfy these conditions with $h=0$.

\subsection{The Main Result}
\label{mainconstruction}
\label{conditions}

We defined $\Psi_{\rm min\/}$ in \S \ref{quantities}.
Given two dyadic objects $Q$ and $\widehat Q$, we define
\begin{equation}
R=\Psi_{\rm min\/}(Q,\widehat Q).
\end{equation}
When $R=0$, our bound gives $\infty$, a result that holds
no matter what.  So, without loss of generality, we treat the
case when $R>0$.    In this case, the two sets
$Q^*$ and $\widehat Q^*$ are contained in disjoint convex sets.

We define
\begin{equation}
\epsilon(Q,\widehat Q)=
\max(0,\Lambda_1)+ \Lambda_2.
\end{equation}
When $Q$ is a dyadic segment, 
\begin{equation}
\Lambda_1=\frac{R}{32} E'(R) + \bigg(\frac{1}{8}-\frac{R^2}{32}\bigg) E''(R);
\hskip 30 pt
\Lambda_2=-\frac{E'(R)}{8}
\end{equation}
When $Q$ is a dyadic square, 
\begin{equation}
\Lambda_1=\frac{R}{16} E'(R) + \bigg(\frac{1}{4}-\frac{R^2}{16}\bigg) E''(R);
\hskip 30 pt
\Lambda_2=-\frac{E'(R)}{7.98} \bigg(\sqrt{1+\overline x^2}+\sqrt{1+\overline y^2}\bigg).
\end{equation}

Let $\delta_i$ be as in Equation \ref{girth}, relative to $Q_i$.  
Recall that $\delta_i$ is a good estimate for the side length of
of the spherical patch $Q_i^*$.
Also, let $\epsilon_{ij}=\epsilon(Q_i,Q_j)$.

\begin{theorem}[Energy Estimator]
\label{EE}
Let ${\cal Q\/}=(Q_0,...,Q_{n})$ be any block.  Then
$${\bf ERR\/}({\cal Q\/})  \leq
\sum_{i=0}^{n-1} \sum_{j \not = i} \epsilon_{ij} \delta_i^2.
$$
\end{theorem}

Again, we remark that the case $n=4$ is the case of interest to us.
For reference, we work out the power-law case $E(r)=r^{-e}$ explicitly.
When $Q$ is a dyadic segment, 
\begin{equation}
\Lambda_1=\frac{e(e+1)}{8R^{e+2}}-\frac{e(e+2)}{32R^e};
\hskip 30 pt
\Lambda_2=\frac{e}{8R^{e+1}}.
\end{equation}
When $Q$ is a dyadic square,
\begin{equation}
\Lambda_1=\frac{e(e+1)}{4R^{e+2}}-\frac{e(e+2)}{16R^e}; \hskip 30 pt
\Lambda_2=\frac{e}{7.98 R^{e+1}} \bigg(\sqrt{1+\overline x^2} +
\sqrt{1+\overline y^2}\bigg)
\end{equation}

\subsection{The Beginning of the Proof}

Recall that a {\it partition of unity\/} is a collection
of functions that sum to $1$ at every point.
A {\it point weighting\/} assigns a partition of unity
\begin{equation}
\lambda_{ab}: Q \to [0,1]; \hskip 30 pt
a,b \in \{0,1\}
\end{equation}
to each dyadic square $Q$ and a partition of unity
\begin{equation}
\lambda_{a}: Q \to [0,1]; \hskip 30 pt
a\in \{0,1\}
\end{equation}
to each dyadic segment.    We have a specific point weighting
in mind, and we will define it in \S \ref{point}.

We let the vertices of a dyadic square $Q$ be
$Q_{ab}$ for $a,b \in \{0,1\}$.  Here
$Q_{00}$ is the lower left vertex and
the remaining vertices, traced in counterclockwise
order, are $Q_{10}$, $Q_{11}$, $Q_{01}$.
We let the vertices of a dyadic segment $Q$ 
be $Q_0$ and $Q_1$, with $Q_0$ on the left.

Given two points $z,w \in \C \cup \infty$, we define
\begin{equation}
f(z,w)=\frac{1}{\|\Sigma^{-1}(z)-\Sigma^{-1}(w)\|^e}.
\end{equation}

Let $\cal X$ denote the set of disjoint pairs
$(Q,\widehat Q)$, as in the previous section.
Let $\epsilon: {\cal X\/} \to \R$ be as in the previous
section.  We will spend
the next $3$ chapters proving the following result.

\begin{lemma}
\label{ABEE}
There exists a point weighting such that the following is true
for all $(Q,\widehat Q) \in \cal X$.

\begin{itemize}

\item When $Q$ is a segment
$$
\bigg(\sum_{a} \lambda_{a}(z) f(Q_{a},w)\bigg)-
f(z,w) < \epsilon(Q,\widehat Q)\ \delta(Q)^2; \hskip 30 pt
\forall (z,w) \in Q \times \widehat Q.
$$

\item When $Q$ is a square
$$
\bigg(\sum_{a,b} \lambda_{ab}(z) f(Q_{ab},w)\bigg)-
f(z,w) < \epsilon(Q,\widehat Q)\ \delta(Q)^2; \hskip 20 pt
\forall (z,w) \in Q \times \widehat Q.
$$
\end{itemize}
\end{lemma}

\noindent
{\bf Proof of Theorem \ref{EE}:\/}
For ease of exposition, we will assume that
the dyadic objects $Q_0,...,Q_{n-1}$ are all dyadic squares.
The case when there are some segments involved 
presents only notational complications.
Consider a configuration
$z_0,...,z_{n}$, domainated by $\cal Q$, that realizes 
${\cal E\/}(\cal Q)$.
For each $i=0,...,n$, define
\begin{equation}
B_i=\{z_0\} \times \{z_{i-1}\} \times Q_i \times ... \times Q_n.
\end{equation}
Just to be explicit, the two extreme cases are
\begin{equation}
B_0=Q_0 \times ... \times Q_n; \hskip 30 pt
B_{n}=\{z_0\} \times ... \times \{z_n\}.
\end{equation}
(Recall that $Q_n=\{z_n\}$.)
The set $B_i$ is a $(2n-2i)$-dimensional 
rectangular solid.    Let ${\cal E\/}_1(i)$ denote the
minimum energy taken over all configurations of $B_i$.
Let ${\cal E\/}_2(i)$ denote the minimum energy over all
vertex configurations of $B_i$.
The conclusion of this lemma is exactly
\begin{equation}
\label{big}
{\cal E\/}_2(0)-{\cal E\/}_1(0) \leq \sum_{i=0}^{n-1}\sum_{j \not = i} 
\epsilon_{ij} \delta(Q_i)^2
\end{equation}
We have
${\cal E\/}_2(n)={\cal E\/}_1(n)$, because $B_{n}$ is a single point.
Hence, we have the telescoping sum
\begin{equation}
\label{induct1}
{\cal E\/}_2(0)-{\cal E\/}_1(0)=
\sum_{i=0}^{n-1}\bigg({\cal E\/}_2(i)-{\cal E\/}_1(i)\bigg) - \bigg({\cal E\/}_2(i+1)-{\cal E\/}_1(i+1)\bigg)
\end{equation}

Our choice of $\{z_i\}$ as a minimizing configuration implies that
${\cal E\/}_1(i)$ is independent of $i$   Combining this with Equation
\ref{induct1}, we get
\begin{equation}
\label{induct3}
{\cal E\/}_2(0)-{\cal E\/}_1(0)=
\sum_{i=0}^{n-1} {\cal E\/}_2(i)-{\cal E\/}_2(i+1).
\end{equation}
To establish Equation \ref{big}, it suffices to prove
\begin{equation}
\label{goal1}
{\cal E\/}_2(i)-{\cal E\/}_2(i+1) \leq \sum_{j \not = i} \epsilon_{ij} \delta_i^2; \hskip 30 pt
\forall i = 0,...,n-1.
\end{equation}

For the remainder of the proof, we fix some
$i \in \{0,...,n-1\}$ once and for all.
We can find vertices $v_{i+1},...,v_{n}$ of $Q_{i+1},...,Q_{n}$
respectively such that 
\begin{equation}
E(z_1,...,z_i,v_{i+1},...,v_{n})={\cal E\/}_2(i+1).
\end{equation}
Here $E(...)$ is the energy of the configuration.
Note that $Q_n$ is a singleton, so that our
``choice'' of $v_n$ is forced.

For any choice of $a,b \in \{0,1\}$, let
$Q_{iab}$ be the corresponding vertex of $Q_i$.
Let \begin{equation}
\lambda_{ab}=\lambda_{ab}(z_i)
\end{equation}
be as guaranteed by Lemma \ref{ABEE}.
The function $\lambda_{ab}$ depends on $i$, but
we suppress this from our notation.
   
Since ${\cal E\/}_2(i)$ is the minimum energy of any vertex
configuration of $B_i$,
\begin{equation}
E(z_1,...,z_{i-1},Q_{iab},v_{i+1},...,v_n) \geq {\cal E\/}_2(i).
\end{equation}
Hence, 
$$
{\cal E\/}_2(i)-{\cal E\/}_2(i+1) \leq A_{ab}; 
$$ 
\begin{equation}
\label{aaa}
A_{ab}=
E(z_1,...,z_{i-1},Q_{iab},v_{i+1},...,v_n)-E(z_1,...,z_i,v_{i+1},...,v_{n}).
\end{equation}

Setting $w_j=z_j$ for $j<i$ and $w_j=v_j$ for $j>i$, we have
\begin{eqnarray}
{\cal E\/}_2(i)-{\cal E\/}_2(i+1) \leq  \cr \cr
\sum_{a,b} \lambda_{ab} A_{ab}=^1 \cr
\sum_{a,b} \lambda_{ab} \bigg(\sum_{j \not = i}  
f(Q_{iab}w_j)-f(z_i,w_j)\bigg)=^2 \cr
\sum_{j \not = i}  \sum_{a,b}
\lambda_{ab} \bigg(f(Q_{iab},w_j)-f(z_i,w_j)\bigg)=^3\cr
\sum_{j \not = i} \bigg(\sum_{a,b} \lambda_{ab} f(Q_{iab},w_j) \bigg)
-f(z_i,w_j)   \leq \cr
\sum_{j \not = i} \epsilon_{ij} \delta_i^2.
\end{eqnarray}

The first inequality comes from Equation \ref{aaa} and
from the fact that $\sum \lambda_{ab}=1$.
Equality 1 comes from the cancellation of all terms not involving
the $i$th index when we subtract the two sums for $A_{ab}$.
Equality 2 comes from switching the order of summation.
Equality 3 comes 
from the fact that $\sum_{ab} \lambda_{ab}=1$.  
The last inequality follows from Lemma \ref{ABEE}.
This establishes Equation \ref{goal1}, which is all we need
to prove Theorem \ref{EE}.
\endproof

\newpage

\section{An Estimate for Line Segments}

\subsection{The Main Estimate}

In this chapter we
prove an estimate that relates directly
to the $\Lambda_1$ term in our Energy
Estimator.  The reader might want to simply
note the main result in this chapter on the
first reading, and then come back to the
proof later on.  Let $E$ be as in the previous chapter.

Let $p \in S^2$ be some point.  In terms
of Lemma \ref{ABEE}, we think of $p$ as being
some point in the spherical patch $(\widehat Q)^*$.  Let
$A' \subset \R^3-\{p\}$ be a segment
whose endpoints lie in $S^2$.  In terms of
Lemma \ref{ABEE}, we think of $A'$ as joining
two boundary points of the spherical patch $Q^*$.
We define
\begin{equation}
\label{energy1}
F(q)=E(\|p-q\|)
\end{equation}
Here $F$ depends on $p$, but we suppress this from our notation.

Define
\begin{equation}
A'_x=(1-x)A'_0+x A'_1; \hskip 30 pt
x \in [0,1].
\end{equation}
Then $x \to A'_x$ is a constant speed
parameterization of $A'$.  

\begin{lemma}[Segment Estimate]
\label{segment}
Suppose that $p$ is at least $R$ units from every point of $A'$.  Let
$\delta$ be the length of $A'$.
$$(1-x) F(A'_0) + x F(A'_1) - F(A'_x)  \leq \frac{X \delta^2}{8},$$
where
$$X=\frac{R}{4} E'(R) + \bigg(1-\frac{R^2}{4}\bigg) E''(R).$$
\end{lemma}

\noindent
{\bf Remark:\/}  We wish to point out one unfortunate
feature of our notation.   The quantities $E'$ and $E''$ are
derivatives of $E$ whereas the quantity $A'$ is simply a 
chord of $S^2$.   We make this notation because, in
the next chapter, we will consider an arc $A$ of $S^2$
and the chord $A'$ that joins the endpoints of $A$.
In some sense, the chord $A'$ is a linear approximation to the
arc $A$
and the derivative $E'$ is a linear approximation to the
function $E$.

\subsection{Strategy of the Proof}

Lemma \ref{segment} really just involves the single variable function
$\phi(x)=F(A'_x)$ defined for $x \in  [0,1]$.    In \S \ref{seg1proof}
we prove the following easy estimate.

\begin{lemma}
\label{seg1}
$$(1-x)\phi(0)+x\phi(1)-\phi(x) \leq \frac{H}{8}; \hskip 30 pt H=\sup_{x \in [0,1]} \phi''(x).$$
\end{lemma}
Here $\phi''(x)$ is the second derivative with respect to $x$.

Recall that $\delta$ is the length of $S$.
Let $s$ denote the arc-length parameter along the segment $S$. 
We set things up so that $s=0$ corresponds to $S_0$ and
$s=\delta$ corresponds to $S_1$.  In general, the parameter $s \in S$
corresponds to $x=s/\delta \in [0,1]$.  
In \S \ref{seg2proof} we establish the following result.

\begin{lemma}
\label{seg2}
Suppose we have the hypotheses in Lemma \ref{segment}.  Then,
$$\frac{d^2 \phi}{ds^2} \leq 
\frac{R}{4} E'(R) + \bigg(1-\frac{R^2}{4}\bigg) E''(R).$$
\end{lemma}

By the Chain Rule, we have
\begin{equation}
\label{seg3}
\phi''(x)=\delta^2\ \frac{d^2 \phi}{ds^2}.
\end{equation}
at corresponding points.  

Lemma \ref{segment} follows from
Lemma \ref{seg1},  Lemma \ref{seg2}, and Equation \ref{seg3}.

\subsection{Proof of Lemma \ref{seg1}}
\label{seg1proof}

We will prove Lemma \ref{seg1} for a smooth function
$h: [0,1] \to \R$.    We want to show that
\begin{equation}
\label{linbound}
(1-x)h(0)+xh(1) - h(x) \leq \frac{H}{8}; \hskip 30 pt H=\sup_{x \in [0,1]} h''(x).
\end{equation}
Let $c$ be any nonzero constant and let $L$ be any linear function.
  Equation \ref{linbound} holds
for the function $h$ if and only if it holds for $ch$.  Likewise,
Equation \ref{linbound} holds for $h-L$ if and only if
Equation \ref{linbound} holds for $h$. Using these two symmetries,
it suffices to prove Equation \ref{linbound} in the case when
$h(0)=h(1)=0$ and $H=1$.
In this case, Equation \ref{linbound} simplifies to
\begin{equation}
- h(x) \leq \frac{1}{8}.
\end{equation}

Let $a \in [0,1]$ be a point where $-h$ attains its maximum.  That
is, $h$ attains its minimum at $a$.  Replacing $h$ by the function
$x \to h(1-x)$ if necessary, we can suppose without loss of generality
$a \geq 1/2$.    

We have $h'(a)=0$ and, by the Fundamental Theorem of Calculus,
\begin{equation}
h'(x) =\int_a^b h''(t)dt \leq x-a
\end{equation}
for all $x \in [a,1]$.
But then
\begin{equation}
-h(a)=h(1)-h(a)=\int_{a}^1 h'(x)dx \leq
\int_a^1 x-a\ dx =
\frac{(1-a)^2}{2} \leq \frac{1}{8}
\end{equation}
This completes the proof.
\endproof

\subsection{Comparison Lemmas}

Let $p$ and $f$ be as above.
Say that a {\it flag\/} is a pair $(q,L)$ where $L$ is a line and
$q \in L$ is a point, and
$q \not = p$.
We define the following quantities.
\begin{itemize}
\item $r(q,L)$ is the distance from $p$ to $q$.
\item $d(q,L)$ is the distance from $p$ to $L$.
\item $\theta(q,L)$ is the small angle between $\overline{pq}$ and $L$.
\item $D(q,L)=F''(q)$, the second derivative w.r.t.
arc length on $L$.
\end{itemize}

We want to compare flats $(q_1,L_1)$ and $(q_2,L_2)$.
We set $\theta_1=\theta(q_1,L_1)$, etc.
Our main result is Corollary \ref{comparison}, proved at the end.
This result will help us establish Lemma \ref{seg2}.

For each $j$ we can rotate so that everything takes place in $\C$, and
$p_j=ir_j$ and $L_j=\R$, and $q_j=x_j$.   Here $x_j^2+d_j^2=r_j^2$.
Suppressing the index, we have
\begin{equation}
\label{partial}
D=\frac{d^2}{dx^2} E(\sqrt{d^2+x^2})\bigg|_x=
E''(r) \frac{x^2}{r^2} + \frac{E'(r)}{r} \times \frac{d^2}{r^2}.
\end{equation}
Here $E'$ and $E''$ are derivatives taken with respect to $r$.

\begin{lemma}
\label{comparison1}
If $\theta_2=\theta_1$ and $r_2 \leq r_1$ and $D_1>0$ then $D_2 \geq D_1$.
\end{lemma}

\startproof
In our situation, there is some constant $\rho \in (0,1]$ so that
\begin{equation}
\frac{r_2}{r_1}=\frac{x_2}{x_2}=\frac{d_2}{d_1}=\rho.
\end{equation}
It follows from Equation \ref{nice} that there are constants
$c \geq 1$ and $h \geq 0$ such that
\begin{equation}
\frac{E'(r_2)}{r_2}=c \frac{E'(r_1)}{r_1}; \hskip 30 pt
E''(r_2)=(c+h) E''(r_1).
\end{equation}
It now follows from Equation \ref{partial} that
\begin{equation}
D_2-cD_1 = h E''(r_1) \geq 0.
\end{equation}
Hence $D_2 \geq c D_1$ for some $c \geq 1$.  
The lemma follows immediately.
\endproof

\begin{lemma}
\label{comparison2}
If $\theta_2 \leq \theta_1$ and $r_2=r_1$ then $D_2 \geq D_1$.
\end{lemma}

The quantity
$d(q,L)$ is monotone increasing with $\theta(q,L)$.  Thus,
we have $d_2 \leq d_1$.     
This time, we have
\begin{equation}
r:=r_2=r_1; \hskip 30 pt
d_2 \leq d_1; \hskip 30 pt x_2 \geq x_1.
\end{equation}
We again have Equation \ref{partial}.
Note that $E''(r)>0$ and $E'(r)<0$.  When we change from
the index $j=1$ to the index $j=2$, we do not decrease the positive
coefficient of $E''(r)$ in the first term and we do not increase the
positive coefficient of $E'(r)/r$ in the second term.
Hence, $D_2 \geq D_1$.
\endproof

The previous two results combine to prove the following corollary.

\begin{corollary}
\label{comparison}
If $\theta_2 \leq \theta_1$ and $r_2 \leq r_1$ and $D_1 \geq 0$, then
$D_2 \geq D_1$.
\end{corollary}

\subsection{Proof of Lemma \ref{seg2}}
\label{seg2proof}

We continue using the notation from above. There is a unique plane
$\Pi$ such that $p \cup L \subset \Pi$.   We rotate the picture so that
$\Pi$ is the $xy$-plane.   Let $C=S^2 \cap \Pi$. Then, as
shown in Figure 6.1, \begin{enumerate}
\item $C$ is a circle
whose radius is at most $1$;
\item  Every point of
$S$ is at least $R$ units from $p$.
\end{enumerate}

\begin{center}
\resizebox{!}{4in}{\includegraphics{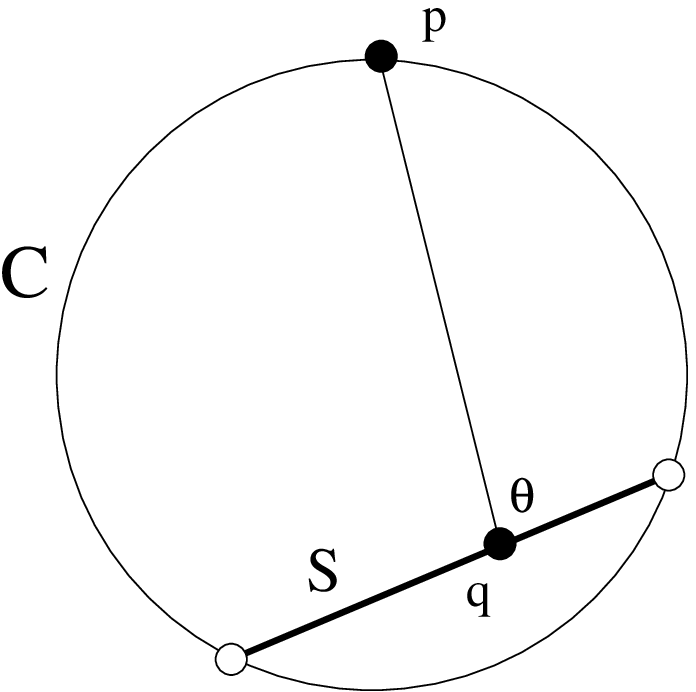}}
\newline
{\bf Figure 6.1:\/} The chord and the circle
\end{center}

We are interested in bounding the quantity $D(q,L)$, where $L$ is the line
containing the chord $S$.   The chord $S$ is subject to the two
constraints mentioned above.  If $C$ has radius less than $1$, we
replace $C$ by a unit radius circle $C'$, and $S$ by a
larger segment $S'$ such that the pair $(C',S')$ satisfies the same
constraints, and the flag $(q,L')$ is the same as the flag $(q,L)$.  Figure
2.2 shows the construction.

\begin{center}
\resizebox{!}{3.7in}{\includegraphics{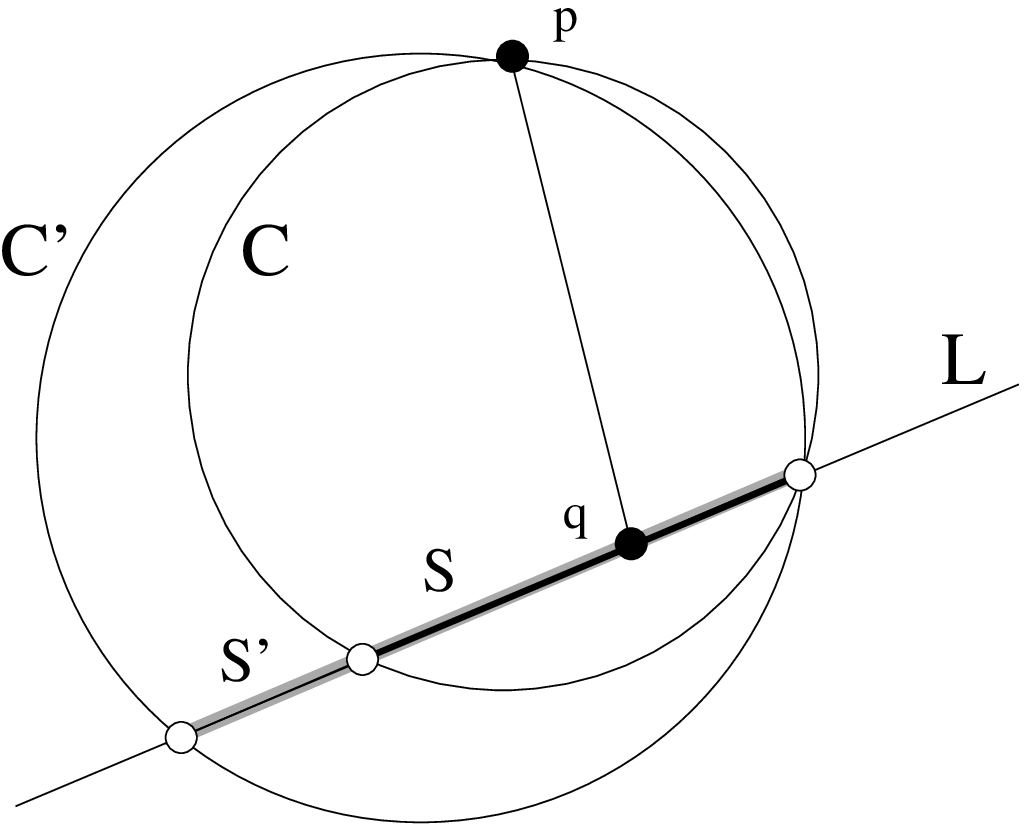}}
\newline
{\bf Figure 6.2:\/} Expanding the circle
\end{center}

So, without loss of generality, we can assume that $C$ is the unit circle.
Also, it suffices to consider the case when $D(q,L) \geq 0$.
Let $q_1=q$.   Let $r_1$ be the distance from $p$ to $q_1$.
Note that $r_1 \geq R$.  Let $q_2$ denote a point on $C$ that
is exactly $R=r_2$ units from $p$.  Let $L_2$ be the line
tangent to $C$ at $q_2$.    We want to apply
Corollary \ref{comparison} to the flats $(q_1,L_1)$ and
$(q_2,L_2)$.  We already know that $r_2 \leq r_1$.

\begin{lemma} $\theta_2 \leq \theta_1$.
\end{lemma}

\startproof
See Figure 6.3.
Let $q_3$ be the endpoint of $S$ such that
the small angle $\theta_1$ subtends the arc
of $C$ between $p$ and $q_3$, as shown in Figure 6.3.
 Let $L_3$ be the line tangent to $C$ at $q_3$.
Let $\theta_3=\theta(q_3,L_3)$.  
The angle $\theta_1$ is half the length
of the two thick arcs in Figure 6.3 whereas the angle
$\theta_3$ is half thelength of the thick arc
joining $p$ to $q_3$.  Hence $\theta_3 \leq \theta_1$.
But the angle $\theta_3$ decreases as we move $q_3$ towards
$p$ along $C$.   Therefore $\theta_2 \leq \theta_3$.
Putting these two inequalities together, we find that
$\theta_2 \leq \theta_1$.
\endproof

\begin{center}
\resizebox{!}{3.8in}{\includegraphics{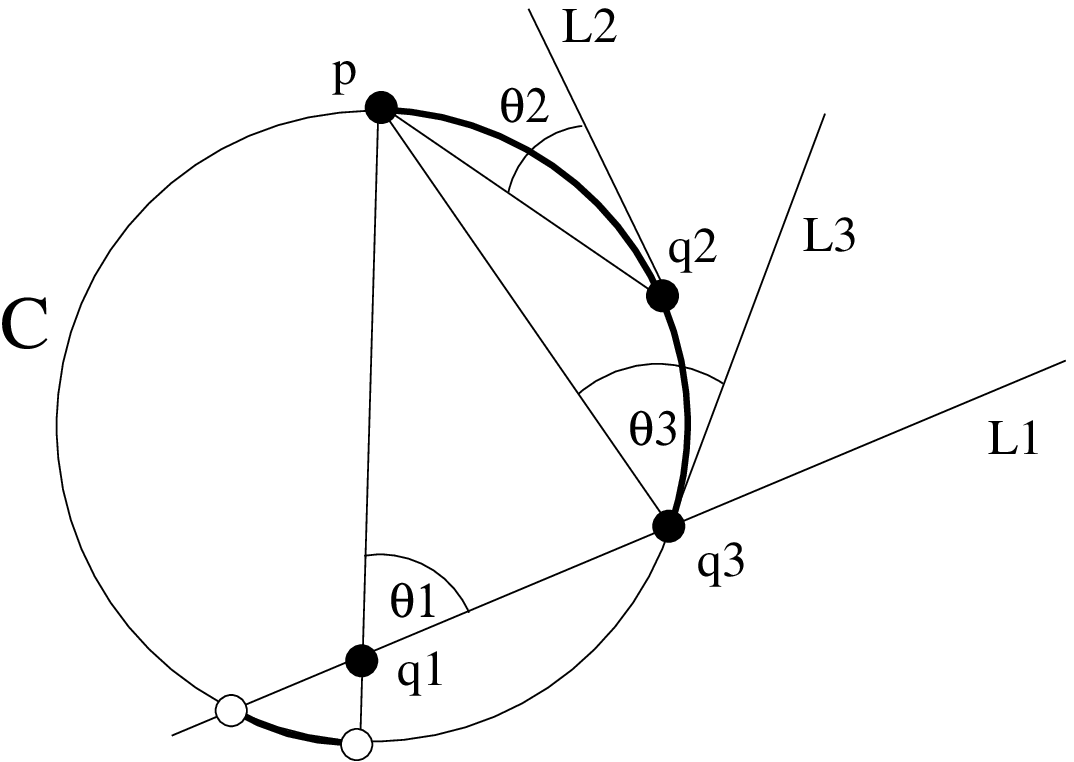}}
\newline
{\bf Figure 6.3:\/} Comparing two flags.
\end{center}

Corollary \ref{comparison} now says that
$D(q_2,L_2) \geq D(q_1,L_1)$.
To finish our proof, it remains only to compute the 
quantity $D(q_2,L_2)$.

We know that $r_2=R$.  Some elementary geometry shows that
\begin{equation}
\label{dval}
d_2=\frac{R^2}{2}; \hskip 30 pt
x_2=r^2-d_2
\end{equation}

Plugging Equation \ref{dval} into Equation \ref{partial} and simplifying,
we get
\begin{equation}
D(q_2,L_2)=
\frac{R}{4} E'(R) + \bigg(1-\frac{R^2}{4}\bigg) E''(R),
\end{equation}
the bound from Lemma \ref{seg2}.  
Now we know that
\begin{equation}
D(q,L) \leq
\frac{R}{4} E'(R) + \bigg(1-\frac{R^2}{4}\bigg) E''(R),
\end{equation}
for all flags $(q,L)$ such that $\|p-q\| \geq R$.  But
$D(q,L)$ is just another name for the quantity
$d^2\phi/ds^2$ featured in Lemma \ref{seg2}.
This completes the proof of Lemma \ref{seg2}.

\newpage

\section{Parametrizing Arcs and Segments}

\subsection{Overview}
\label{overview}

Let $S \subset \C$ be a line segment.   In this chapter we
define a certain parametrization of $S$ by a parameter
$x \in [0,1]$.  
Once we define our parametrization, we will state and prove
several geometric results about it.  As with the last chapter,
the reader might want to just note the results on the first
reading and then come back later for the proofs.

Let $A$ be a circular arc on $S^2$, and let $A'$ be the
chord that joins the endpoints $A_0$ and $A_1$ of $A$.
Throughout the chapter, we assume that $A$ is contained
in a semicircle.
We let $x \to A'_x$ be the affine map from $[0,1]$ to $A'$.
Let $C$ be the circle containing $A$, and let $c \in C$ be the
point which is diametrically opposed to the midpoint of $A$.
We define $A_x$ so that the three points
$c,A'_x,A_x$ are always collinear.  Here is our first result.

\begin{lemma}
\label{circle2}
The function $f(x)=\|A_x-A'_x\|$ attatains its maximum at $x=1/2$.
\end{lemma}

Now we return to our main task of parametrizing a segment
$S \subset \C$.  Let
$S^*=\Sigma^{-1}(S)$.  The method above
gives us a parametrization of $S^*$.  Now we define
\begin{equation}
\label{parameter}
S_x=\Sigma(S^*_x).
\end{equation}
As usual $\Sigma$ denotes stereographic projection.

Suppose $Q$ is a normal dyadic square.
This means that $Q$ does not cross the coordinate axes,
and the side length of $Q$ is at most $1$.
We consider the case when $Q$ is contained in the positive
quadrant.  The other cases have symmetric treatments.
Let $Q_0$ and $Q_1$ be the left
and right edges of $Q$.   For any $x \in [0,1]$,
let $Q_x$ denote the segment connecting
$(Q_0)_x$ and $(Q_1)_x$.    The main result in this
chapter gives estimates on the size and shape of
the image of $Q_x^*=\Sigma^{-1}(Q_x)$.

\begin{lemma}
\label{circle3}
$Q_x^*$ has arc length at most 
$$\delta \times 1.0013.$$
and is contained in a circle of radius at least
$$\frac{1}{\sqrt{1+\overline y^2}}.$$
\end{lemma}

\subsection{Proof of Lemma \ref{circle2}}

Our result is scale-invariant.   It suffices to prove the
result when $A$ is an arc of the unit circle, as shown
in Figure 7.1.  The arc $cy$ is evidently shorter than
the diameter $cx$.    On the other hand, the arc
$cz$ is evidently longer than the arc $cw$.  Hence
the arc $yz$ is shorter than the arc $wx$.
This is what we wanted to prove.

\begin{center}
\resizebox{!}{3in}{\includegraphics{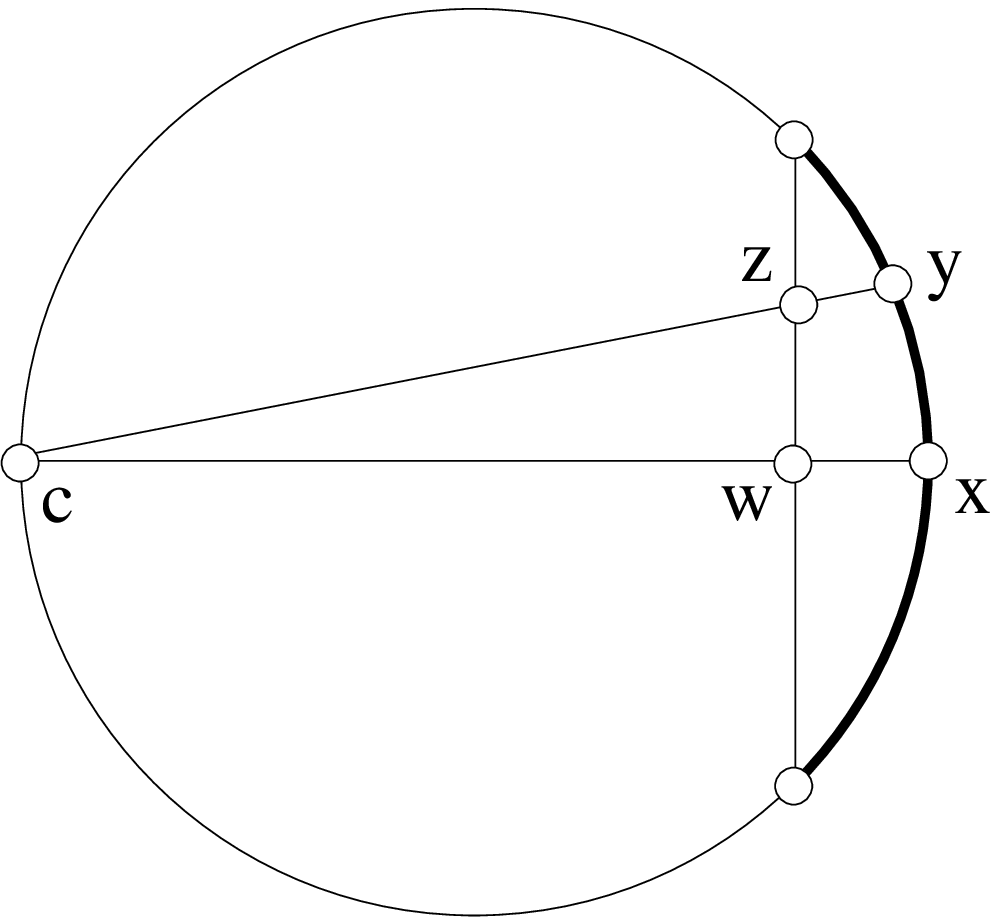}}
\newline
{\bf Figure 7.1:\/} The relevant points
\end{center}

\subsection{A Lemma about Slopes}

The rest of the chapter is devoted to proving Lemma \ref{circle3}.
Here we reduce Lemma \ref{circle3} to a more geometric statement.

\begin{lemma}
\label{circle4}
$Q_x$ has slope in $(0,0.051)$ for all $Q$ and $x \in (0,1)$.
\end{lemma}

\noindent
{\bf Proof of Lemma \ref{circle3}:\/}
Let $S=Q_x$.
We deal first with the arc length of $S^*$.
This is really the same argument as in Lemma \ref{confine}.
When evaluated at all points of $S$, the quantity in
Equation \ref{metric2} is at most $\delta/s$.
Since $S$ has slope in $(0,0.051)$ and $Q$ has side length $s$,
the segment $S$ has length at most 
\begin{equation}
s \times \sqrt{1+(0.051)^2}<1.0013.
\end{equation}
The arc length estimate on $S^*$ follows by integration.

Since the line $L$ through $S$ has positive slope, it intersects
the imaginary axis in a point of the form $iy$ where
$y<\underline y$.  But then, according to Equation \ref{metric1},
there are two points on $(L \cup \infty)^*$ which are at least
$$\frac{1}{\sqrt{1+\overline y^2}}$$
apart.
\endproof

The rest of the chapter is devoted to proving Lemma \ref{circle4}.

\subsection{Mobius Geometry}
\label{mob}
Let the vertices of $Q$ be
$Q_{00},Q_{10}, Q_{11}, Q_{01}$, starting in the bottom
left corner and going counterclockwise around.
Associated to $Q$ is an auxilliary map
$\phi: [0,1] \to [0,1]$, defined as follows.
\begin{itemize}
\item The point $(1-x)Q_{00} + x Q_{01}$ is $(Q_0)_s$ for some $s$.
\item The point $(Q_1)_s$ is $(1-y) Q_{10}+y Q_{11}$ for some $y=\phi(x)$.
\end{itemize}
In other words, we consider the natural map $(Q_0)_s \to (Q_1)_s$ but
we precompose and postcompose with affine maps to make the
domain and range equal to $[0,1]$.  
Lemma \ref{circle3} is equivalent to the statement that
\begin{equation}
0 <\phi(x) -x <0.051; \hskip 30 pt \forall x \in (0,1).
\end{equation}

Given $4$ distinct points $A,B,C,D \in \R^n$, we have
the cross ratio
\begin{equation}
\chi(A,B,C,D)=\frac{\|A-C\|\ \|B-D\|}{\|A-B\|\ \|C-D\|}.
\end{equation}
We say that a homeomorphism from one curve to another is
{\it Mobius\/} if the map preserves cross ratios.
In case the curves are line segments in the plane, a Mobius
map between them is the restriction of a linear fractional
transformation of $\C \cup \infty$. 

\begin{lemma}
$\phi$ is Mobius.
\end{lemma}

\startproof
Since similarities are Mobius transformations, it suffices to prove
that the map $(Q_0)_s \to (Q_1)_s$ is a Mobius transformation.
Stereographic projection is well known to be a Mobius map
from any line segment in $\C$ to the corresponding arc on $S^2$.
Referring to the construction in the beginning of
\S \ref{overview}, the map from $A_x$ to $A'_x$ is just the
composition of affine maps with (one dimensional) stereographic
projection.  Hence, this map is also Mobius.  

Let $A_0=Q_0^*$ and $A'_0$ be the chord connecting the
endpoints of $A_0$.  Likewise define $A_1$ and $A'_1$.  The
map of interest to us is the composition
\begin{equation}
\label{bigcomp}
Q_0 \to A_0 \to A'_0 \to A'_1 \to A_1 \to Q_1.
\end{equation}
The outer maps are Mobius, from what we have already said,
and the middle map is affine.
\endproof

A Mobius map from $[0,1]$ to $[0,1]$ is completely determined
by its derivative at either endpoint.  Thus, we can understand
$\phi$ by computing or estimating $\phi'(0)$.   
In our next result, we think of $\phi'(0)$ as a function of the
choice of dyadic square $Q$.     As in Equation \ref{metric2},
define
\begin{equation}
g(x,y)=\frac{2}{1+x^2+y^2}.
\end{equation}
Define
\begin{equation}
G_Q=\frac{g(Q_{00})g(Q_{11})}{g(Q_{10})g(Q_{10})}.
\end{equation}

\begin{lemma}
Relative to $Q$, we have
$$\phi'(0)=\sqrt{G_Q}.$$
In particular, $\phi$ is the identity iff $G_Q=1$.
\end{lemma}

\startproof
If follows from symmetry that
\begin{equation}
\phi'(0)\phi'(1)=1.
\end{equation}
Therefore
\begin{equation}
\label{root1}
\phi'(0)=\sqrt{\frac{\phi'(0)}{\phi'(1)}}.
\end{equation}
Looking at the composition in
Equation \ref{bigcomp}, all the maps except the
outer two have the same derivative at either endpoint.
For the affine map in the middle, this is obvious: the
derivative is constant.
In the case of the map $A_j \to A'_j$ this follows from the fact that we
are projecting from a point $c_j$ that is symmetrically located
with respect to $A'_j$ and $A_j$.  Call this property of
the derivatives the {\it symmetry property\/}.

By Equation \ref{metric2}, the quantity $g(Q_{00})$ is the
norm of the derivative of $\Sigma^{-1}$ at $Q_{00}$.
The other quantities $g(Q_{ij})$ have similar interpretations.
It therefore follows from the symmetry property and the
Chain Rule that
\begin{equation}
\label{root2}
\frac{\phi'(0)}{\phi'(1)}=\frac{g(Q_{00})/g(Q_{01})}{g(Q_{10})/g(Q_{11})}.
\end{equation}
This Lemma now follows from Equations \ref{root1} and \ref{root2}.
\endproof

\subsection{The End of the Proof}

For the purposes of doing calculus, we define $\phi$ and $G$ relative
to any square that is contained in the positive quadrant. 
 We remind the reader that we only consider squares that have
side length at most $1$.

\begin{lemma}
It never happens that $G_Q=1$.
\end{lemma}

\startproof
When $Q_{00}=(x,y)$ and $Q$ has side length $r$, we compute
\begin{equation}
\label{mama}
G_Q=\frac{(1+r^2+2rx+x^2+y^2)(1+r^2+2ry+x^2+y^2)}
{(1+x^2+y^2)(1+2r^2+x^2+y^2+2rx+2ry)}
\end{equation}
Every factor in Equation \ref{mama} is a polynomial with only
constant coefficients and at least one constant term.  Hence
this expression never vanishes.
\endproof

When $Q=[0,1]^2$, we compute that $G_Q=4/3$.   Hence
$\phi'(0)>1$ relative to this choice of $Q$.  It tollows
from continuity that $\phi'(0)>1$ relative to any
square in the positive quadrant.
But this means that $\phi$ is increasing.  (Here, of
course, we are crucially using the fact that $\phi$ is
a Mobius map from $[0,1]$ to $[0,1]$.) Hence
$\phi(x)>x$.   This proves that the slope of the
segment $Q_x$ is positive for all $Q$ and all $x>0$.
This proves half of Lemma \ref{circle4}.
Now we turn to the other half.

\begin{lemma}
\label{calculus}
The quantity $G_Q$ is maximized when $Q$ has side length $1$ and
$$Q_{00}=(\xi,\xi); \hskip 30 pt \xi=\frac{\sqrt 3-1}{2}.$$
\end{lemma}

\startproof
Let $\psi(x,y,r)$ be the function in Equation \ref{mama}.
We compute symbolically that 
\begin{equation}
\frac{d\psi}{dr}=\frac{\Delta(x,y,r)}{(1+x^2+y^2)^2(1+2r^2+x^2+y^2+2rx+2ry)^2},
\end{equation}
Where $\Delta$ is a polynomial with entirely positive terms, at least one of which
involves only $r$.  (The polynomial is rather long and unenlightening.)
From this we conclude that $d\psi/dr>0$.   Hence, the maximum value
of $G_Q$ must occur when $r=1$.

We compute that $\psi(\xi,\xi,1)=3/2$.   So, we just have to show that
the function $h(x,y)=3/2-\psi(x,y,1)$ is non-negative on the positive quadrant.
We compute that $h$ is a rational function.  The denominator is a polynomial
with only positive terms, and the numerator $N$ equals
\begin{equation}
1-2x+4x^2+2x^3+x^4-2y-8xy+2x^2y+4y^2+2xy^2
+2x^2y^2+2y^3+y^4
\end{equation}
In fact, $N$ is non-negative on the entire plane.  To see this,
we make the change of variables
\begin{equation}
x=\xi+u; \hskip 50 pt y=\xi+v.
\end{equation}
With this change of variables, we find that
$N-(2u-2v)^2$ is a polynomial involving only
positive terms.
\endproof

In light of the previous result, the quantity $\phi'(0)$ is maximized
for the special square $Q_0$ from Lemma \ref{calculus}.
Since $G=3/2$ in this case, we have
\begin{equation}
\phi'(0)=\sqrt{3/2}.
\end{equation}
But this equation pins down $\phi$ uniquely, and we observe that the map
\begin{equation}
\phi(x)=\frac{3 \sqrt 2 x}{2 \sqrt 3 + 3 \sqrt x - 2 \sqrt 3 x}.
\end{equation}
has the same derivative.  Hence, this is the correct formula for $\phi$.
A bit of calculus now shows that
\begin{equation}
\phi(x)-x<0.051; \hskip 30 pt \forall x \in [0,1].
\end{equation}
This completes the proof.

\newpage

\section{Proof of Lemma \ref{ABEE}}
\label{point}

\subsection{The Geometry of Circles}

We need one more result about circles.

\begin{lemma}
\label{circle}
Suppose that $A$ is an arc of a circle $C$.  Let
$d$ be the arc-length of $A$.  Let
$A'$ be the segment connecting the endpoints of $A$.
Let $r$ be the radius of $C$.     Let
$\mu$ be the smallest constant 
such that every point of $A$ is within
$\mu$ of some point of $A'$.  Then
$$\mu< \frac{d^2}{8r}.$$
\end{lemma}

\startproof
Let's first consider the case $r=1$.
We rotate so that $C$ is the unit circle, and $A$ is the
arc bounded by the points $\exp(-i\theta)$ and $\exp(i \theta)$.
Here $\theta \in (0,\pi)$.   Then 
\begin{equation}
d=2 \theta; \hskip 30 pt
\mu=1-\cos(\theta).
\end{equation}
The claim of this lemma boils down to the statement that
\begin{equation}
\frac{\theta^2}{1-\cos(\theta)}>2,
\end{equation}
This is equivalent to the statement that
\begin{equation}
\phi(\theta)=\theta^2+2\cos(\theta)-2>0.
\end{equation}
We have $\phi(0)=0$ and
\begin{equation}
\phi'(\theta)=2(\theta-\sin(\theta)>0.
\end{equation}
This proves what we need.

If $r \not = 1$, we let $T$ be a dilation that scales distances
by a factor of $1/r$.   The arc $T(A)$ has length $d/r$ and
$T(C)$ has radius $1$. We apply our result to the pair
$(T(A),T(C))$ and find that every point of $T(A)$ is within
$$\frac{(d/r)^2}{8}$$ 
of $T(C)$.   Applying $T^{-1}$, we get the desired result for the
pair $(A,C)$.
\endproof

\subsection{Lemma \ref{ABEE} for dyadic segments}

\subsubsection{Defining the Weighting}

Suppose that $(Q,\widehat Q)$ is a reasonable pair,
and $Q$ is a dyadic line segment.  The construction in the
previous chapter gives us a parametrization $x \to Q_x$.
We take $Q_0$ to be the left
endpoint and $Q_1$ to be the right endpoint.  
The corresponding endpoints of $Q^*$ are
$Q_0^*$ and $Q_1^*$.
We define
our weighting as follows.  Letting $z=Q_x$, we define
\begin{equation}
\lambda_0(z)=1-x; \hskip 30 pt
\lambda_1(z)=x;
\end{equation}

\subsubsection{Setting up the Calculation}

To bring our notation in line with Lemma \ref{segment}, we define
\begin{equation}
A=A(Q_0^*,Q_1^*)=Q^*; \hskip 30 pt
A'=A'(Q_0^*,Q_1^*).
\end{equation}

Let $(z,w) \in Q \times \widehat Q$.   Define
\begin{equation}
p=\Sigma^{-1}(w); \hskip 30 pt q=\Sigma^{-1}(z).
\end{equation}
Setting $x=\lambda_1(z)$, we have
\begin{equation}
q=A_x.
\end{equation}
We also define
\begin{equation}
q'=A'_x
\end{equation}
The idea of our proof is to estimate things with $q'$ in place of $q$, and
then to estimate the error we get when replacing $q'$ by $q$.

\subsubsection{Using Lemma \ref{segment}}

Let $F$ be as in Equation \ref{energy1}.
We have
$$
\sum_u \lambda_u(z) f(Q_u,w)=(1-x)F(A'_0)+x F(A'_1);$$
\begin{equation}
F(A'_x)=E(\|q'-p\|).
\end{equation}
Lemma \ref{segment} now tells us that
\begin{equation}
\label{proof0}
\sum_u \lambda_u(z) f(Q_u,w) -E(\|p-q'\|)  \leq \max(0,\Lambda_1)\ \delta^2
\end{equation}

\subsubsection{The Easy Case}

Now we need to see what happens when we replace $q'$ by $q$.
Define
\begin{equation}
 r=\|p-q\| \hskip 30 pt r'=\|p-q'\|.
\end{equation}
If $r' \geq r$ then (since $E$ is decreasing)
\begin{equation}
\label{easy} 
E(\|p-q\|) \geq E(\|p-q'\|).\end{equation}
In this case,
our proof is done: Equations \ref{easy} and \ref{proof0} 
combine to give a tighter bound than what Lemma \ref{ABEE} gives.

\subsubsection{The Hard Case}

Now suppose $r'<r$.
Since $E$ is convex, $E'$ is monotone decreasing.  Therefore
\begin{equation}
\label{proof1}
E(r')-E(r) \leq E'(r') \|q-q'\| \leq E'(R)\|q-q'\|.
\end{equation}
Here $R \geq r'$ is as in Lemma \ref{ABEE}.

The same proof as in Lemma \ref{circle3} shows that $A$ has arc length
at most $\delta$.   Also $A$ is contained in a great circle -- i.e. a circle
of radius $1$.
Lemma \ref{circle} now says that every point of $A$ is within
$\delta^2/8$ of $A'$.    But the point of $A'$ closest to
$A_{1/2}$ is $A'_{1/2}$.  Therefore
\begin{equation}
\label{proof2}
\|q-q'\|=\|A_x-A'_x\| \leq^* \|A_{1/2}-A'_{1/2}\| \leq \frac{\delta^2}{8}.
\end{equation}
The starred inequality is Lemma \ref{circle2}.
Combining Equations \ref{proof1} and \ref{proof2}, we find that
\begin{equation}
\label{proof3}
E(r')-E(r) \leq E'(R) \frac{\delta^2}{8} = \Lambda_2\  \delta^2.
\end{equation}
Note that
\begin{equation}
\label{proof4}
f(z,w)=E(\|p-q\|).
\end{equation}
Hence, by Equation \ref{proof4}, 
\begin{equation}
\label{proof5}
E(\|p-q'\|)-f(z,w)=E(r')-E(r)\leq \Lambda_2\ \delta^2.
\end{equation}
Adding Equations \ref{proof0} and \ref{proof5}, we get
the bound in Lemma \ref{ABEE}.    This completes the proof
in case $Q$ is a dyadic segment.

\subsection{The Weighting for Dyadic Squares}

Let $Q$ be a dyadic square and let $Q^*=\Sigma^{-1}(Q)$.
Let $Q_{ab}$ be the vertices of $Q$, as in
\S \ref{mob}.  Let $Q^*_{ab}$ be the corresponding vertex of $Q^*$.
As we make our construction, the reader should picture the
letter `H', with the horizontal bar very
slightly slanted.  We will use coordinates $(h,v)\in [0,1]^2$.
The $v$ variable moves along vertical segments and the
$h$ variable moves along (roughly) horizontal segments.

Let $Q_0$ be the left edge of $Q$.   The two endpoints of
$Q_0$ are $Q_{00}$ and $Q_{01}$.  Likewise, let
$Q_1$ be the right edge of $Q$.   The two endpoints of
$Q_{10}$ are $Q_{11}$.   Using the parametrization
from the previous chapter, we define
\begin{equation}
Q_{0v}=(Q_0)_v\in Q_0; \hskip 30 pt
Q_{1v}=(Q_1)_v\in Q_1.
\end{equation}
Next, we define  $Q_v$ to be the segment joining $Q_{0v}$ to $Q_{1v}$.
Note that $Q_v$ is nearly horizontal but not exactly horizontal.
See Lemma \ref{circle4}.
Finally, we define
\begin{equation}
Q_{hv}=(Q_v)_h,
\end{equation}
again using the parametrization discussed in the previous sections.
We define our weighting as follows.   Assuming that $z=Q_{vh}$,
$$
\lambda_{00}(z)=(1-h)(1-v); \hskip 30 pt
\lambda_{10}(z)=(h)(1-v); 
$$
\begin{equation}
\lambda_{01}(z)=(1-h)(t); \hskip 30 pt
\lambda_{11}(z)=(h)(v).
\end{equation}

In the previous section, we defined a weighting for
dyadic segments.  We could equally well define this
weighting for any segment, since we have already
explained how to parametrize any such segment.
There are three ways in which our weighting here is compatible with
the weighting for segments.

\begin{itemize}
\item With respect to the segment $Q_0$, the weighting for the point
$S_{0v}$ is $\lambda_0=1-v$ and $\lambda_1=v$.
\item With respect to the segment $Q_1$, the weighting for the point
$S_{1v}$ is $\lambda_0=1-v$ and $\lambda_1=v$.
\item With respect to the segment $Q_v$, the weighting for the
point $Q_{hv}$ is $\lambda_0=1-h$ and $\lambda_1=h$.
\end{itemize}

\subsection{The End of The Proof}

For ease of notation, we will assume that our dyadic square
lies in the positive quadrant.  The other cases are similar,
and indeed follow from symmetry.
We gather $4$ pieces of information.
\begin{enumerate}

\item The circles containing the arcs $Q_0^*$ and $Q_1^*$ 
have radius at least
$$\frac{1}{\sqrt{1 + \overline x}}.$$
This follows from Equation \ref{metric1} and from the fact that
the line extending a vertical edge of $Q$ comes within $\overline x$ of the origin.

\item The arcs $Q_0^*$ and $Q_1^*$ have length at most
$\delta$.   This follows from the same argument as in Lemma \ref{circle3}.

\item For any $v\in [0,1]$, the line
extending the segment $Q_v$ comes within
$\overline y$ of the origin.  Hence, the circle
containing $Q_v^*$ has radius at least
$$\frac{1}{\sqrt{1 + \overline y}}.$$
See Lemma \ref{circle3}.

\item For any $v \in [0,1]$, the arc $Q_v^*$ has length
at most $(1.0013) \times \delta$.  See Lemma \ref{circle3}.

\end{enumerate}

Now we are ready for the main argument.  The basic idea is to make
repeated appeals to the segment case of Lemma \ref{ABEE} and
then to suitably average the result.

We define
$$
\Lambda_{1x}=\Lambda_{1y}=\Lambda_1/2;
$$
\begin{equation}
\Lambda_{2x}=-\frac{E'(R)}{8} \sqrt{1+\overline x^2}; \hskip 30 pt
\Lambda_{2y}=-\frac{E'(R)}{7.98} \sqrt{1+\overline y^2}.
\end{equation}
We have
\begin{equation}
\Lambda_1=\Lambda_{1x}+\Lambda_{1y}; \hskip 30 pt
\Lambda_2=\Lambda_{2x}+\Lambda_{2y}.
\end{equation}

Let $p \in (\widehat Q)^*$ be some point.
Let $w=\Sigma(p)$.
Let $f$ be as in Lemma \ref{ABEE}.  That is
\begin{equation}
f(z,w)=E(\|\Sigma^{-1}(z)-\Sigma^{-1}(w)\|).
\end{equation}

Applying the segment case of Lemma \ref{ABEE} to the arc
$Q_0$, we find that
\begin{equation}
\label{add1}
Z_1:=(1-v)f(Q_{00},w)+vf(Q_{01},w)-f(Q_{0v},w) \leq
\bigg(\max \Lambda_{1x}+\Lambda_{2x}\bigg) \delta^2
\end{equation}
The proof is exactly the same as in the previous section, except for the
one point that the circular arc $Q_0^*$ lies not necessarily in a
great circle but rather a circle whose radius is bounded by Item 1 above.
The designation $Z_1$ is for algebraic purposes which will become
clear momentarily.
Similarly
\begin{equation}
\label{add2}
Z_2:= (1-v)f(Q_{10},w)+vf(Q_{11},w)-f(Q_{1v},w) \leq
\bigg(\max \Lambda_{1x}+\Lambda_{2x}\bigg) \delta^2
\end{equation}

Finally, an argument just like the one given for {\it dyadic\/} segments
also works for the segment $Q_s$.  The only property we used about
dyadic segments is that they don't cross the coordinate axes, and
$Q_s$ has this property.    Using Items 3 and 4 above in place of
Items 1 and 2, the same argument gives

\begin{equation}
\label{add3}
Z_3:= (1-h)f(Q_{0v},w)+hf(Q_{1v},w)-f(Q_{hv},w) \leq
\bigg(\max \Lambda_{1y}+\Lambda_{2y}\bigg) \delta^2
\end{equation}
The key point here is that 
$$7.98<8/1.0013.$$
The number $1.0013$ comes up in Item 4 above.

Concentrating on the left hand sides of
Equations \ref{add1}, \ref{add2}, and \ref{add3}, we have
$$(1-h)Z_1+hZ_2 + Z_3 =$$
$$ (1-h)(1-v)f(Q_{00},w)+
 (1-h)(v)f(Q_{01},w)+$$
$$ (h)(1-v)f(Q_{10},w)+
(h)(v)f(Q_{11},w)-
f(Q_{hv},w)=$$
\begin{equation}
\label{quant}
\bigg(\sum_{ab} \lambda_{ab}(z) f(Q_{ab},w)\bigg)- f(z,w); \hskip 30 pt
z=Q_{hv}.
\end{equation}
Concentrating on the right hand sides, we see that the
quantity in Equation \ref{quant} is at most
$$ (1-h) \bigg(\max \Lambda_{1x}+\Lambda_{2x}\bigg) \delta^2+
 (h) \bigg(\max \Lambda_{1x}+\Lambda_{2x}\bigg) \delta^2+
\bigg(\max \Lambda_{1y}+\Lambda_{2y}\bigg) \delta^2=$$
\begin{equation}
\bigg(\max(\Lambda_1,0)+\Lambda_2\bigg) \delta^2.
\end{equation}
This proves Lemma \ref{ABEE}.

\section{The Hessian and its Variation}
\label{hessian}
\label{derproof}

\subsection{Main Part of the Proof}

In this chapter we prove Lemma \ref{mainhess}.   
Let $H_e$ denote the Hessian of $\cal E$, the energy function,
relative to the function $E(r)=r^{-e}$.
Let $Z$ be the configuration corresponding to the TBP, normalized
as in Equation \ref{min1}. 

\begin{lemma}
\label{posdef}
The lowest eigenvalue of $H_e$ exceeds $1/10$ for $e=1,2$.
\end{lemma}

\startproof
Let $M$ be either of the two matrices.  Let
$I_7$ be the $7 \times 7$ identity matrix. Using a modified version of
the Cholesky Decomposition, as discussed in [{\bf Wa\/}, p 84], we
write
\begin{equation}
(M-\frac{1}{10} I_7)=L D L^t.
\end{equation}
where $L$ is lower triangular, $D$ is diagonal matrix with all
diagonal entries positive, and $L^t$ is the transpose of $L$.
This suffices to show that $M-\frac{1}{10} I_7$ is positive
definite.  Hence, the lowest eigenvalue of $M$ is at least $1/10$.
\endproof

Given a square matrix $M$, we define
\begin{equation}
\|M\|_2 =\bigg( \sum_{i,j} M_{ij}^2\bigg)^{1/2}.
\end{equation}
We mention a familiar and useful property of this norm.
\begin{lemma}
\label{cs}
For any unit vector $v \in \R^7$ we have
$\|M(v) \cdot v\| \leq \|M\|_2$.
\end{lemma}

\startproof
We have
$$
\|M(v) \cdot v\|_2=
\bigg{\|}\sum M_{ij} v_iv_j \bigg{\|}_2 \leq^*
\|M\|_2 \bigg{\|}\sum(v_iv_j)^2\bigg{\|}^{1/2}=$$
\begin{equation}
\|M\|_2\sqrt{\bigg(\sum v_i^2\bigg)\bigg(\sum v_j^2\bigg)}=
\|M\|_2 \|v\|^2=\|M\|_2.
\end{equation}
The starred inequality is the Cauchy-Schwarz inequality.
\endproof

We write
$M' \prec M$ if $|M'_{ij}| \leq M_{ij}$ for all indices.
Recall that $s=2^{-11}$, as in Lemmas \ref{maincomp}
and \ref{mainhess}.
Let $D_k$ denote the partial derivative with respect
to the $k$th direction in $\R^7$.   

\begin{lemma}[Variation Bound]
\label{derbound}
Let $\Psi_{1,e},...,\Psi_{7,e}$ be the smallest non-negative matrices such that
$D_kH_e(Z)  \prec \Psi_{k,e}$ for all $k$ and all $Z \in \Omega_{s}$.
Let
$$\Psi(e)=\bigg{\|}\sum_{k=1}^7 \Psi_{k,e}\bigg{\|}_2.$$
Then $\Psi(1)<345$ and $\Psi(2)<140$ and, 
$\sup_e \Psi(e)<463$.
\end{lemma}

\begin{corollary}
Let $Z$ be the center of $\Omega_s$ (i.e., the TBP)
 and let $W \in \Omega_s$ be
any other point. Then
$$\|H_e(Z)-H_e(W)\|_2<\frac{1}{10}; \hskip 30pt
\forall e \in (0,\infty).$$
\end{corollary}

\startproof
Let $H=H_e$.
Say that a {\it special path\/} in $\Omega$ is a $7$-segment
polygonal path $\gamma$ such that the $i$th segment is
parallel to the $i$th coordinate direction.
We can connect $Z=Z_0$ to $W=Z_7$ by a special path
$\gamma$, all of whose segments have length at
most $2^{-12}$.
Let $Z_0,...,Z_7$ be the vertices of $\gamma$.
Let 
$\Delta_{ij}=H(Z_i)-H(Z_j)$.  Integrating along the $k$th segment,
we get the bound
$\Delta_{k,k-1} \prec  2^{-12}\Psi_k$.  Hence
\begin{equation}
\|H(Z)-H(W)\|_2=\|\Delta_{70}\|_2=
\bigg{\|}\sum_{k=1}^7 \Delta_{k,k-1} \bigg{\|}_2 
\leq  \frac{\bigg{\|}\sum_{k=1}^7 \Psi_k \bigg{\|}_2}{2^{12}}<
\frac{345}{4096}.
\end{equation}
This last quantity is less than $1/10$.
\endproof

\noindent
{\bf Proof of Lemma \ref{mainhess}:\/}
Let $e=1$ or $e=2$.  Let $H=H_e$.
Let $W \in \Omega_s$ be some point.
We want to show that
$H(W)$ is positive definite.
Let $A=H(Z)$ and $B=H(W)$.
By the preceding Corollary, we have $B=A+\Delta$
with $\|\Delta\|_2<1/10$.

Now let $v \in \R^7$ be any unit vector.
We have
\begin{equation}
Bv\cdot v = Av\cdot v + \Delta v \cdot v \geq
Av \cdot v - |\Delta v \cdot v|>\frac{1}{10}-\|\Delta\|_2>0.
\end{equation}
Hence $B$ is positive definite.
The fact that $Av \cdot v>1/10$ comes from the fact that
the lowest eigenvalue of the symmetric matrix $A$ is greater than $1/10$.
\endproof

\subsection{A Stronger Bound}

Consider the function
\begin{equation}
\label{powerf}
\phi(x)=x^{e/2}
\end{equation}
and the interval
\begin{equation}
\label{interval0}
I=\left[\frac{1}{4}+\frac{1}{2}+2^{-9}\right].
\end{equation}
Define
\begin{equation}
\label{defn}
c_k(e)=\sup_{x \in I} \frac{d^k\phi}{dx_k}(x)
\end{equation}
Setting $c_k=c_k(e)$, we define
\begin{equation}
\label{ups}
\Upsilon(e)=
\sqrt{19336 c_1^2+19036 c_1c_2 + 4922 c_2^2+ 1474 c_1c_3 + 772 c_2c_3 + 31 c_3^2}.
\end{equation}

Our next result uses the notation from the Variation Bound.
\begin{lemma}[Variation Bound II]
\label{derbound2}
$\Psi(e)<\Upsilon(e)$ for all $e \in (0,\infty)$.
\end{lemma}

We find easily that
\begin{equation}
(c_1(1),c_2(1),c_3(1))=(1,2,12).
\end{equation}
Plugging this into Equation \ref{ups}, we get $\Psi(1)<345$. 
Similarly, we have
\begin{equation}
(c_1(2),c_2(2),c_3(2))=(1,0,0).
\end{equation}
This yields $\Psi(2)<140$.

Some elementary calculus shows that
\begin{equation}
\sup_{e \in (0,\infty)} c_1(e)<1.1; \hskip 30 pt
\sup _e c_2(e)<3.2; \hskip 30 pt
\sup_{e \in (0,\infty)} c_3(e)<16.
\end{equation}
We omit the details, because we don't use the general bounds anywhere
in our main proof.  When we plug in these bounds we find that
$\Psi(e)<463$ for all $e$.  Hence, the Variation Bound II implies
the Variation Bound.
\newline
\newline
{\bf Remark:\/}
The Variation Bound II is generally much better than the
Variation Bound.    For instance,
$\Upsilon(e) \to 0$ as $e \to 0$ or as $e \to \infty$.

\subsection{Proof of the Variation Lemma II}

Our proof usually suppresses the dependence on the exponent $e$.

 A point $(z_0,z_1,z_2,z_3) \in \Omega_s$ is such that
each $z_m$ lies within a square $\Delta_{m+1}$ of side length $2^{-11}$ about
one of the points of the TBP configuration, normalized as
in Equation \ref{min1}.  Let $R_1$ be the union of these
$4$ squares, $\Delta_1,...,\Delta_4$.
Let $R_2 \subset \C^2$ denote the set of points
$z_1,z_2$ which arise as a disjoint pair of finite points of
a configuration of $\Omega_{s}$.   Here $R_2$ consists
of $12$ components, all of the form
$\Delta_i \times \Delta_j$.  The
components $\Delta_i \times \Delta_j$
and $\Delta_j \times \Delta_i$ are {\it partners\/}.
We let $\Delta_5,...,\Delta_{10}$ be $6$ components,
no two of which are partners.

Recall that $\Sigma^{-1}$ is inverse stereographic projection.
See Equation \ref{INV}.  Setting $z=x+iy$, define
\begin{equation}
\label{FF}
F(z)=\frac{1}{\|\Sigma^{-1}(z)-\Sigma^{-1}(\infty)\|^2}=
\frac{1+x^2+y^2}{4}.
\end{equation}

\begin{equation}
G(z_1,z_2)=
\frac{1}{\|\Sigma^{-1}(z_1)-\Sigma^{-1}(z_2)\|^2}=
\frac{(1+x_1^2+y_1^2)(1+x_2^2+y^2)}{4(x_1-x_2)^2+(y_1-y_2)^2}.
\end{equation}

Let $\phi$ be as in Equation \ref{powerf}.  According
as $m \leq 4$ or $m \geq 5$, define
\begin{equation}
\widehat U_m=\phi \circ U_m; \hskip 30 pt
U_m=F
\hskip 5 pt {\rm or\/} \hskip 5 pt G
\end{equation}
Define
\begin{equation}
\label{components}
\Phi(i,j,k,m)=\sup_{\Delta_m} |D_iD_jD_k U_m|; \hskip 20 pt
\widehat \Phi(i,j,k,m)=\sup_{\Delta_m} |D_iD_jD_k \widehat U_m|.
\end{equation}
Here $i,j,k \in \{1,...,7\}$ and $m \in \{1,...,10\}$ and
$D_k$ is the $k$th partial derivative.

We have
\begin{equation}
{\cal E\/}(z_0,z_1,z_2,z_3)=\sum_{m=1}^{10} \widehat U_m,
\end{equation}
where the arguments of the functions on the right hand side
are suitably chosen sub-lists of $(z_0,z_1,z_2,z_3)$.  Therefore
\begin{equation}
|D_k D_iD_j({\cal E\/})| \leq  \Psi_k(i,j):= \sum_{m=1}^{10} \widehat \Phi(i,j,k,m).
\end{equation}
With this definition of $\Psi_k$, 
we have $D_kH \prec \Psi_k$.
Summing over $k$, we have
\begin{equation}
\label{vb2}
\bigg{\|}\sum_{k=1}^7 \Psi_k\bigg{\|}_2 \leq
\bigg{\|}\sum_{k,m} \widehat \Phi(i,j,k,m)\bigg{\|}_2.
\end{equation}

Say that $i \in \{1,...,7\}$ is {\it related\/} to an index $m \in \{1,....,10\}$ if
a change in the coordinate $x_i$ moves one of the points in the argument
of $\widehat U_m$.  For instance $i=2$ is related to $m=5$ because changing
$x_2$ changes the location of the point $z_1=x_1+ix_2$, and the argument
of $U_5$ is the point $(z_0,z_1) \in \Delta_5$.   Clearly
$\widehat \Phi(i,j,k,m)=0$ unless $i,j,k$ are all related to $m$.

When all indices are related to $m$, we use the chain rule
$$
D_iD_jD_k \widehat U_m =
 \phi' \times D_iD_jD_kU_m+\phi''' \times D_i U_m \times D_j U_m \times D_k U_m+$$
\begin{equation}
\label{cr0}
\phi'' \times D_iD_j U_m \times  D_kU_m +
\phi'' \times D_jD_k U_m  \times  D_iU_m +
\phi'' \times D_kD_i U_m  \times  D_j U_m.
\end{equation}

Given the definitions of the regions $R_1$ and $R_2$ we have
$U_m(\Delta_m) \subset I$,
where $I$ is the interval from Equation \ref{interval0}.
Combining this fact with
Equation \ref{defn} and Equation \ref{cr0}, we have

$$
\widehat \Phi(i,j,k,m) \leq
c_1 \Phi(i,j,k,m) + c_3 \Phi(i,m) \Phi(j,m)\Phi(k,m)+$$
\begin{equation}
\label{cr1}
 c_2 \Phi(i,j,m)\Phi(k,m)+
 c_2 \Phi(j,k,m)\Phi(i,m)+
 c_2 \Phi(k,i,m)\Phi(j,m).
\end{equation}

Let $\delta_m$ be the center of $\Delta_m$.   Define
\begin{equation}
a(i,j,k,m)=|D_iD_jD_kU_m(\delta_m)|.
\end{equation}
We define $a(i,j;m)$ and $a(i;m)$ similarly. Define
\begin{equation}
v=\big((v(i),v(i,j),v(i,j,k)\big)=\bigg(\frac{1}{1000},0,0\bigg) \hskip 15 pt
{\rm or\/} \hskip 15 pt
\bigg(\frac{1}{200},\frac{1}{50},\frac{1}{10}\bigg),
\end{equation}
according as $m \leq 4$ or $m \geq 5$.  
Let $b(i;m)=a(i;m)+v(i)$, etc.

\begin{lemma}
\label{estimate}
for all $(i,j,k,m)$ we have
$$\Phi(i,m) \leq b(i,m); \hskip 20 pt
\Phi(i,j,m) \leq b(i,j,m); \hskip 20 pt
\Phi(i,j,k,m) \leq b(i,j,k,m).$$
\end{lemma}

Combining Equation \ref{vb2}, Equation \ref{cr1}, and
Lemma \ref{estimate}, we have
\begin{equation}
\sum_k \Psi_k  \prec
\sum_{k,m} \delta(i,j,k,m)
\left(\matrix{c_1  b(i,j,k,m) + \cr \cr
c_2  (b(i,j,m)  b(k,m) + \cr\cr
c_2  (b(j,k,m)  b(i,m) + \cr\cr
c_2  (b(k,i,m)  b(j,m) + \cr\cr
c_3  b(i,m)  b(j,m)  b(k,m) }\right)
\end{equation}
Here $\delta(i,j,k,m)=1$ if $i,j,k$ are all related to $m$,
and otherwise $0$.
When we compute the norm of the right hand side
in Mathematica, we get the square root of a polynomial
in $c_1,c_2,c_3$.   Rounding the coefficients of this
polynomial up to integers, we get $\Upsilon$.
This completes the proof of the Variation Bound II,
modulo the proof of Lemma \ref{estimate}.

\subsection{Proof of Lemma \ref{estimate}}

When $m \leq 4$, we are dealing with the function $F$
from Equation \ref{FF}.  
We have $D_i F=x_i/2$ and all higher derivatives are constant.
Lemma \ref{estimate} in this trivial case now follows from the fact that
${\rm radius\/}(\Delta_m)<1/500$.

For the nontrivial cases,
we fix some value of $m \in \{5,...,10\}$ and consider
the cube $\Delta=\Delta_m \subset \C^2$.
This cube has side-length $2^{-11}$.
The center point is $\delta=\delta_m$.  Define
\begin{equation}
a_k=|\D_k G(\delta)|; \hskip 30 pt
\Phi_k=\max_{(z_1,z_2) \in \Delta} |\D_k G(z_1,z_2)|.
\end{equation}
The maxima here are taken over all $k$th partial derivatives $\D_k$.
An exact, finite calculation, done for each of the
$6$ choices of $m$, yields
\begin{equation}
\label{explicit}
a_1 \leq 1; \hskip 30 pt
a_2 \leq 4; \hskip 30 pt
a_3 \leq 18; \hskip 30 pt
a_4 \leq 96; \hskip 30 pt
a_5 \leq 600.
\end{equation}

Now we give an estimate on $\Phi_6$.
Given a polynomial $P \in \R[x_1,y_1,x_2,y_2]$, let
$|P|$ denote the sum of the absolute values of the coefficients of $P$.
We have the easy upper bound
\begin{equation}
|P(z_1,z_2)| \leq |P| \max(|z_1|,|z_2|)^d; \hskip 30 pt
d={\rm degree\/}(P).
\end{equation}
Here we have set $z_j=x_j+iy_j$.
Calculating symbolically we find that each $6$th partial derivative $D_6$ has
the following structure.  There is a polynomial $P$, depending on the
choice of derivative, such that
\begin{equation}
\label{form}
D_6G = \frac{P}{\|z_1-z_2\|^7}; \hskip 30 pt
|P| \leq 4519440; \hskip 30 pt
{\rm deg\/}(P) \leq 10.
\end{equation}
We (easily) 
have
\begin{equation}
\max(\|z_1\|,\|z_2\|, \|z_1-z_2\|^{-1})<1+2^{-8}; \hskip 30 pt
\forall (z_1,z_2) \in \Delta.
\end{equation}
Hence
\begin{equation}
\Phi_6 \leq 4519440  \times (1+2^{-8})^{17}<5000000.
\end{equation}

\begin{lemma}
\label{bootstrap}
$\Phi_k<a_k+2^{-10} \times \Phi_{k+1}$.
\end{lemma}

\startproof
Let $(w_1,w_2)$ be a point of $\Delta$.   We can connect
$(w_1,w_2)$ to a point $(z_1,z_2) \in \delta_2$ by a
$4$-segment polygonal path $\gamma$ such that
the $i$th segment has length at most $2^{-12}$ and
moves in the $i$th coordinate direction.
This lemma now follows from integration.
\endproof

Now we apply Lemma \ref{bootstrap} in an iterative way.

\begin{equation}
\Phi_5 \leq 600 + 5000000 \times  2^{-10}<5483
\end{equation}

\begin{equation}
\Phi_4 \leq 96 + 5483 \times 2^{-10} <102.
\end{equation}

\begin{equation}
\label{l1}
|\Phi_3-a_3| < 102 \times 2^{-10}<\frac{1}{10}; \hskip 30 pt \Phi_3<18.1
\end{equation}

\begin{equation}
\label{l2}
|\Phi_2-a_2| < 18.1 \times 2^{-10}<\frac{1}{50}; \hskip 30 pt
\Phi_2<4.02.
\end{equation}

\begin{equation}
|\Phi_1-a_1| \leq  4.02 \times 2^{-10}<\frac{1}{200}.
\end{equation}
The R.H.S. of Equation \ref{l1} comes from the L.H.S. and Equation \ref{explicit}.
Likewise, the R.H.S. of Equation \ref{l2} comes from the L.H.S. and Equation \ref{explicit}.
This completes the proof of Lemma \ref{estimate}.

\newpage

\section{Computational Issues}
\label{interval}

\subsection{General Remarks}

Our calculations in this paper are of two kinds.
The material in \S \ref{hessian} makes some exact calculations in
Mathematica [{\bf W\/}] and the rest of the paper makes calculations in Java.
For the Mathematica calculations, we need to take derivatives of
rational functions or their square roots, and evaluate them at
elements of $\Q[\sqrt 3]$.   We also need to simplify and
group terms of some polynomials.   Everything is manipulated
exactly.  The user who downloads our Java program will also
find our Mathematica files in the same directory.

The bulk of the calculations are done in Java, and these
are what we discuss below.
We take the IEEE-754 standards for binary floating point arithmetic
[{\bf I\/}] as our reference for the Java calculations.  This 1985 document
has recently been superceded by a 2008
publication [{\bf I2\/}].   We will stick to the 1985 publication for
three reasons.  First, the 1985 version is shorter and simpler.
Second, the portion relevant to our computation has not changed in
any significant way.  Third, we think that
some of the computers running our code will have been manufactured
between 1985 and 2008, thus conforming to the older standard.

\subsection{Doubles}

With a view towards explaining interval arithmetic, we
first describe the way that Java represents real numbers.
Our Java code represents real numbers by {\it doubles\/}, in a way that is
an insignificant modification of the scheme
discussed in [{\bf I\/}, \S 3.2.2].     To see what our program does,
read the documentation for the
{\bf longBitsToDouble\/} method in the {\bf Double\/} class,
on the website
{\bf http://java.sun.j2se/1.4.2.docs/api\/}.
According to this documentation -- and experiments verify that it works
this way on our computer -- our program represents a double
by a $64$ bit binary string, where
\begin{itemize}
\item The first bit is called $s$.
\item The next $11$ bits are a binary expansion of an integer $e$.
\item If $e=0$, the last $52$ bits are the binary expansion of an integer $m$.
\item If $e \not = 0$, the last $52$ bits are binary expansion of the number
$m-2^{52}$.
\end{itemize}
The real number represented by the double is
\begin{equation}
(-1)^s  \times 2^{e-1075} \times m.
\end{equation}

\noindent
{\bf Example:\/}
The double representing $-317$ is
stored as the $64$ bit string
$$1/10000000111/001111010 \ldots 0.$$
The slashes are put in to emphasize the breaks.
Here $s=1$ and (since $e \not = 0$)
$$e=2^{10}+8+4+1=1031.$$
There are $44$ zeros at the end of the word, and
$$m=2^{44}(0+0+32+16+8+4+0+1)+2^{52}=
2^{44}(317).$$
Hence 
$$(-1)^s \times 2^{e-1075} \times m = -317.$$

Now we come to the main point of our discussion above.
The non-negative doubles have a lexicographic ordering, and this ordering
coincides with the usual ordering of the real numbers they
represent.   The lexicographic ordering for the non-positive doubles
is the reverse of the usual ordering of the real numbers they
represent.  To {\it increment\/} $x_+$ of
 a positive double $x$ is the very next double
in the ordering.  This amounts to treating the last $63$ bits of the
string as an integer (written in binary) and adding $1$ to it.
With this interpretation, we have $x_+=x+1$.
We also have the decrement $x_-=x-1$.
Similar operations are defined on the non-positive doubles.
These operations are not defined on the largest and smallest
doubles, but our program never encounters (or comes anywhere near)
 these.

\subsection{The Basic Operations}
\label{basicop}

Let $\D$ be the set of all doubles.
Let 
\begin{equation}
\R_0=\{x \in \R|\ |x| \leq 2^{30}\}
\end{equation}
Our choice of $2^{30}$ is an arbitrary but convenient cutoff.
Let $\D_0$ denote the set of doubles representing reals in $\R_0$.

According to the discussion in 
[{\bf I\/}, 3.2.2, 4.1, 5.6], there is a map
$\R_0 \to \D_0$ which maps each $x \in \R_0$ to some
$[x] \in \D_0$ which is closest to $x$.   In case there
are several equally close choices, the computer chooses one
according to the method
in [{\bf I\/}, \S 4.1].    This ``nearest point projection''
exists on a subset of $\R$ that is much larger
than $\R_0$, but we only need to consider
$\R_0$.  We also have the inclusion $r: \D_0 \to \R_0$, which
maps a double to the real that it represents.   

Our calculations use the $5$ functions
\begin{equation}
+; \hskip 20 pt
-; \hskip 20 pt
\times; \hskip 20 pt
\div; \hskip 20 pt
\sqrt{\hskip 5 pt}.
\end{equation}
These operations act on $\R_0$ in the usual way.
Operations with the same name act on $\D_0$.
Regarding the first $5$ basic operations, [{\bf I\/}, \S 5] states that
{\it each of the operations shall be performed as if it first produced an
intermediate result correct to infinite precision and with unbounded range, and then 
coerced this intermediate result to fit into the destination's format\/}.  Thus,
for doubles $x$ and $y$.
\begin{equation}
\label{rule}
\sqrt{x}=\bigg[\sqrt{r(|x|)}\bigg]; \hskip 20 pt
x*y=[r(x)*r(y)]; \hskip 20 pt
* \in \{+,-,\times,\div\}.
\end{equation}
The operations on the left hand side represent operations on doubles
and the operations on the right hand side represent operations on
reals.
\newline
\newline
{\bf Remark:\/}
Exceptions to Equation \ref{rule} 
can arise if we divide by $0$ or a number
too close to zero.  This will produce an {\it overflow error\/} -- too
large a number to be accurately represented by a double.
To see that this never happens, we check that we never divide
by a number smaller than $2^{-11}$.     Moreover,
as we discuss at the end of this chapter,
we never take the square root of a number less than $2^{-5}$.
\newline

\subsection{Interval Arithmetic}

An {\it interval\/} is a pair $I=(x,y)$ of doubles with $x \leq y$.   Say
that $I$ {\it bounds\/} $z \in \R_0$ if $x \leq [z] \leq y$.  This is true
if and only if $x \leq z \leq y$.  Define
\begin{equation}
[x,y]_o=[x_-,y_+].
\end{equation}
This operationis well defined for doubles in $\D_0$.
We are essentially {\it rounding out\/} the endpoints of the interval.
Let $I_0$ and $I_1$ denote the left and right endpoints of $I$.
Letting $I$ and $J$ be intervals, we define
\begin{equation}
\label{operate}
I*J = (\min_{ij} I_i *I_j,\max_{ij} I_i*I_j)_o; \hskip 30 pt
\sqrt I=(\sqrt I_0,\sqrt I_1)_0
\end{equation}
That is, we perform the operations on all the endpoints, order
the results, and then round outward. 
Given Equation \ref{rule}, we have the following facts.
\begin{enumerate}
\item If $I$ bounds $x$ and $J$ bounds $y$ then
$I*J$ bounds $x*y$.
\item If $I>0$ bounds $x>0$ then $\sqrt I$ bounds $\sqrt x$.
\end{enumerate}

After we have defined intervals and their basic operations, we
define other Java objects based on intervals.  Namely,
interval versions of complex numbers, vectors in $\R^3$, 
dyadic segments and squares, and dyadic boxes.  For instance,
the interval version of a complex number is an object of the form
$X+iY$, where $X$ and $Y$ are intervals.
The algebra of these interval objects -- e.g. complex addition or 
the dot product -- is formally the
same as the corresponding algebra on the usual objects.  At every
step of the calculation, the real version of the object is bounded,
component by component, by the interval version.  If some
particular interval object passes our test, it means that all the
real objects bounded by it also pass.

\subsection{Implementation Details}

We implement the interval arithmetic in a way that
tries to minimize the time we spend using it.
We run the floating point algorithm until we notice that
a box has passed one of the floating point tests.
Then, we re-perform the interval arithmetic test on the
interval arithmetic version of the box.  
\begin{itemize}
\item  If the box passes
the interval arithmetic test, we eliminate it, and switch back
to the floating point algorithm.  
\item If the box fails the interval arithmetic test, we just act as if it
failed the corresponding floating point test, and resume the algorithm.
We call this case a {\it mismatch\/}.  It is harmless.
\end{itemize}

Even though the mismatches are harmless, we prefer to
have few mismatches, so that run time for the interval arithmetic
version of the code is easier to predict from the run time of the
floating point version.   There are $4$ kinds of potential
mismatches, corresponding to the $4$ kinds of tests we
perform, as discussed in \S \ref{discuss1}.   The two
mismatches that actually arise (or, rather arose) in
abundance are
the Energy Estimator mismatches
and the Redundancy Eliminator mismatches.
Once we make our modifications, the code runs
completely without mismatches for the $1/r$ potential
and with only a $4$ mismatches for the
$1/r^2$ potential.

To prevent the Energy Eliminator mismatches,
we arrange the floating point eliminator so that a dyadic
box passes only if the minimum energy is above
$E_0+2^{-40}$, where $E_0$ is the energy of the T.B.P.
In this way, the interval versions of the boxes that are
actually checked have a bit of a cushion.
The reason why this fudge factors does not
cause our code to halt is that the inequality
${\cal E\/}(X)-E_0<2^{-40}$ only occurs well
inside the $L^{\infty}$ neighborhood of size
$2^{-14}$ about the TBP.

Adding the fudge factor of $2^{-40}$ makes
the floating point test harder to pass, but
 does not change the
validity of the program.   Logically, the floating
point calculation simply finds a candidate
partition of the configuration space, in which
each box in the partition passes one of our tests.
Adding the fudge factor does change the final
partition a bit, but it doesn't destroy the basic
property of the partition.

Mismatches occur for the Redundancy Eliminator
because we sometimes eliminate configurations
where there is an exact equality of coordinates.
This equality will fail for the corresponding
intervals, because of a tiny overlap.
To get around this problem, we associate to each dyadic object a
Gaussian integer that represents $2^{25} z$, where $z$
is the center of the object.   When we subdivide a dyadic
object, we perform the arithmetic on the Gaussian integer
exactly, so as to compute the exact value of the centers
of the subdivided objects.  We never subdivide more than
$24$ times -- and in fact the maximum is about $17$ -- so
we never arrive at a situation where $2^{25} z$ is not 
an integer.  
When it comes time to compare the various
coordinates of a dyadic object, we actually compare
$2^{25}$ times those coordinates, so that we are
working entirely with integers.  This lets us make
exact comparisons.
\newline

There is one more fine point we would like to mention.
We want to avoid taking expressions
of the form $\sqrt I$, where $I$ is an interval that is too close to
$0$.   If $0 \in I$, then $I_0<0$ and $\sqrt I_0$ causes an
arithmetic error.  The only time this issue comes up is in the
bounds from Lemma \ref{bound}.  The quantity $4-D^2$ might
be close to $0$ if $D$ is close to $2$.  To avoid this problem,
we define a new function $\sigma$, which has the property
$\sigma(I)=2^{-5}$ if $I_1<2^{-10}$ and otherwise
$\sigma(I)=\sqrt I$.  (It never happens that $I_1>2^{-10}$ and
$I_0<0$.) When computing the interval version
of the bound in Lemma \ref{bound}, we use $\sigma$ in
place of $\sqrt{\hskip 5 pt}$.  This causes no problem with the proof,
because the bound is not as strong when we use
$\sigma$ in place of $\sqrt{\hskip 5 pt}$.

\newpage

\section{References}

\noindent
[{\bf BBCGKS\/}] Brandon Ballinger, Grigoriy Blekherman, Henry Cohn, Noah Giansiracusa, Elizabeth Kelly, Achill Schurmann, \newline
{\it Experimental Study of Energy-Minimizing Point Configurations on Spheres\/}, 
arXiv: math/0611451v3, 7 Oct 2008
\newline
\newline
[{\bf C\/}] Harvey Cohn, {\it Stability Configurations of Electrons on a Sphere\/},
Mathematical Tables and Other Aids to Computation, Vol 10, No 55,
July 1956, pp 117-120.
\newline
\newline
[{\bf CK\/}] Henry Cohn and Abhinav Kumar, {\it Universally 
Optimal Distributions of Points on Spheres\/}, J.A.M.S. {\bf 20\/} (2007) 99-147
\newline
\newline
[{\bf CCD\/}] online website: \newline
http://www-wales.ch.cam.ac.uk/$\sim$ wales/CCD/Thomson/table.html
\newline
\newline
[{\bf DLT\/}] P. D. Dragnev, D. A. Legg, and D. W. Townsend,
{\it Discrete Logarithmic Energy on the Sphere\/}, Pacific Journal of Mathematics,
Volume 207, Number 2 (2002) pp 345--357
\newline
\newline
[{\bf HS\/}], Xiaorong Hou and Junwei Shao,
{\it Spherical Distribution of 5 Points with Maximal Distance Sum\/}, 
arXiv:0906.0937v1 [cs.DM] 4 Jun 2009
\newline
\newline
[{\bf I\/}] IEEE Standard for Binary Floating-Point Arithmetic
(IEEE Std 754-1985)
Institute of Electrical and Electronics Engineers, July 26, 1985
\newline
\newline
[{\bf I2\/}] IEEE Standard for Floating-Point Arithmetic
(IEEE Std 754-2008)
Institute of Electrical and Electronics Engineers, August 29, 2008.
\newline
\newline
[{\bf RSZ\/}] E. A. Rakhmanoff, E. B. Saff, and Y. M. Zhou,
{\it Electrons on the Sphere\/},  Computational Methods and Function Theory,
R. M. Ali, St. Ruscheweyh, and E. B. Saff, Eds. (1995) pp 111-127
\newline
\newline
[{\bf SK\/}] E. B. Saff and A. B. J. Kuijlaars,
{\it Distributing many points on a Sphere\/}, 
Math. Intelligencer, Volume 19, Number 1, December 1997 pp 5-11
\newpage
\noindent
[{\bf T\/}] J. J. Thomson, {\it On the Structure of the Atom: an Investigation of the
Stability of the Periods of Oscillation of a number of Corpuscles arranged at equal intervals around the
Circumference of a Circle with Application of the results to the Theory of Atomic Structure\/}.
Philosophical magazine, Series 6, Volume 7, Number 39, pp 237-265, March 1904.
\newline
\newline
[{\bf W\/}] S. Wolfram, {\it The Mathematica Book\/}, 4th ed. Wolfram Media/Cambridge
University Press, Champaign/Cambridge (1999)
\newline
\newline
[{\bf Wa\/}] D. Watkins, {\it Fundamentals of Matrix Calculations\/},  John Wiley and Sons (2002)
\newline
\newline
[{\bf Y\/}], V. A. Yudin, {\it Minimum potential energy of a point system of charges\/}
(Russian) Diskret. Mat. {\bf 4\/} (1992), 115-121, translation in Discrete Math Appl. {\bf 3\/} (1993) 75-81

\end{document}